%% file: semconvex_arxiv.tex
\documentclass[a4paper]{article}
\usepackage[text={7in,10.3in},centering]{geometry}
\usepackage{graphicx} 
\usepackage{epstopdf}
\usepackage{subcaption}
\usepackage{float}
\usepackage{color}
\usepackage{multicol,multirow}
\usepackage{listings}
\usepackage{setspace}
\usepackage{lineno}

\usepackage{hyperref}

\usepackage{amsmath}
\usepackage{amsthm}
\usepackage{amssymb}
\usepackage{gensymb}	 

\newtheorem{theorem}{Theorem}
\newtheorem{lemma}[theorem]{Lemma}
\newtheorem{proposition}[theorem]{Proposition}

\graphicspath{{figures/}} 
\input zdefs

\begin{document}
\title{Convex Formulation for Regularized Estimation \\ of Structural Equation Models}

\author{Anupon Pruttiakaravanich and Jitkomut Songsiri\footnote{Corresponding author}  \\  Department of Electrical Engineering, Faculty of Engineering \\ Chulalongkorn University \\ \texttt{anupon106@gmail.com, jitkomut.s@chula.ac.th}}
\maketitle

\begin{abstract}
Path analysis is a model class of structural equation modeling (SEM), which it describes causal relations among measured variables in the form of a multiple linear regression. This paper presents two estimation formulations, one each for confirmatory and exploratory SEM, where a zero pattern of the estimated path coefficient matrix can explain a causality structure of the variables. The original nonlinear equality constraints of the model parameters were relaxed to an inequality, allowing the transformation of the original problem into a convex framework. A regularized estimation formulation was then proposed for exploratory SEM using an l1-type penalty of the path coefficient matrix. Under a condition on problem parameters, our optimal solution is low rank and provides a useful solution to the original problem. Proximal algorithms were applied to solve our convex programs in a large-scale setting. The performance of this approach was demonstrated in both simulated and real data sets, and in comparison with an existing method. When applied to two real application results (learning causality among climate variables in Thailand and examining connectivity differences in autism patients using fMRI time series from ABIDE data sets) the findings could explain known relationships among environmental variables and discern known and new brain connectivity differences, respectively.

\end{abstract}

\newpage
\section{Introduction}
\label{sec:intro}
Structural equation modelling (SEM) is a class of multivariate models used for learning a causal relationship among variables (exploratory modelling) or for testing whether the model is best fit by given data (confirmatory modelling). A general SEM includes the observed and latent variables, while their relationships are explained by a linear model whose parameters explain the cause or influence from one variable to another. As such, SEM has been widely used in behavioural research, such as in psychology, sociology, business and medical research~\cite{McL:94, PLFI:09}. The details of the model and its historical background in SEM are given in~\cite[\S 1]{Bol:89}. Path analysis is a special problem in SEM, where it provides a model for explaining the relationships among measured variables only (no latent variables). This can be better associated with scientific research where the observed variables are often of primary interest. For example, one aims to explore causal relationships among brain regions from brain signals (such as fMRI data)~\cite{McL:94,BuF:97,BHHB+:00,KZCB+:07,JKCH+:09,CGSH+:11} where the entries of the path coefficient matrix in the model explain how much change in activities of one region influences another region. 

In path analysis, one applies a prior knowledge about the relationship structure of variables of interest to construct a model and encode a structure, such as the zero pattern of the \emph{path matrix}. The first problem type in path analysis, the~\emph{confirmatory SEM}, is to estimate the value of nonzero entries in the path matrix and the covariance matrix of the model residual errors so that the model-reproduced covariance matrix fits well with the sample covariance matrix in an optimal sense, as evaluated by various types of criterion functions, such as maximum likelihood (ML) and ordinary or weighted least-squares (LS)~\cite[\S 4]{Bol:89}. The second type of problem,~\emph{exploratory SEM}, is to learn the causal structure of variables from a zero structure of the estimated path matrix. The existing approach for exploratory SEM is to begin with a base model, where a certain set of paths are affirmative but the existence of some other paths is in question. This results in a set of a few candidate models associated with different zero structures of the path matrix and then the significance of the difference between these models can be determined from the $\chi^2$ statistic~\cite[\S 7]{Bol:89}. Examples of this approach can be seen in brain network studies~\cite{McL:94,BuF:97}, where only a few variables (in the order of up to $10$ brain regions) are selected. One can locally search for a path structure by starting from a null model and sequentially allowing the coefficient corresponding to the largest Lagrangian multiplier to be nonzero~\cite{BHHB+:00}. The most optimal, but not tractable, method is to perform an exhaustive search that enumerates all possible path patterns and chooses the model corresponding to the lowest minimized ML function~\cite{CGSH+:11}. However, as the number of all possible models grows exponentially to the number of variables, it is not a feasible approach. 

Both confirmatory and exploratory SEM problems are nonlinear optimization problems in matrix variables with a quadratic equality and positive definite cone constraints. Common techniques based on Newton-Raphson or gradient descent are implemented to estimate the model parameters~\cite[\S 7]{Mul:09}, \cite[\S 4]{Bol:89} and there are many existing SEM commercial softwares, such as LISREL, EQS and Mplus~\cite{RaM:06,JSTT:00,Bol:89}, so in an estimation process, a starting value for the update iteration is required. Although these numerical methods work well under normal conditions, some initial values may not lead to the convergence of the optimal solution or may stick into local minima, and so several strategies for selecting initial values have been proposed~\cite[\S 4]{Bol:89}. These include choosing an instrumental variable estimate or selecting the strength of the path coefficient magnitude. When the iterative method in these software does not converge, the user is suggested not to interpret the result. 

In this work, we present two alternative estimation formulations, one each for confirmatory and exploratory SEM problems. The original nonlinear equality constraints of the model parameters are relaxed to an inequality, leading the problems to transform into convex formulations that can be solved efficiently by many existing convex program solvers, where the solution is guaranteed to be the global minima. For exploratory SEM, we propose an objective function that is added with an $\ell_1$-type regularization of the path coefficient matrix, called sparse SEM. Such formulations are regarded as a \emph{lasso} formulation~\cite{HTF:09} and encourage many zeros in the path matrix, allowing us to read off the zero pattern and interpret it as a causal structure of the variables. Solving exploratory SEM using a regularization approach has previously been discussed~\cite{JGM16}, where the reticular action model (RAM) was the focus of the model class, or where the model explains a relationship between the dependent variables and latent variables~\cite{ZTT:16, HCW:17}.

Difficulties in estimating sparse coefficients in latent equations arise~\cite{ZTT:16, HCW:17}, and so the expected conditional maximization (ECM) algorithm was applied to manage the unobserved latent variables. To our knowledge, the algorithms in~\cite{ZTT:16,HCW:17} may potentially suffer from local minimum issues if one attempts to solve high-dimensional settings, such as in brain applications. Despite the fact that these model classes are more general than a path analysis problem, our approach could serve as a convex framework targeted to solve a special class of RAM, and was shown to perform better than~\cite{JGM16} in some cases. We applied alternating direction method of multipliers (ADMM), and proximal parallel algorithm (PPXA), which requires a feasible amount of memory storage suitable for large-scale implementation. More importantly, we showed that, under a condition on problem parameters, the optimal covariance error is diagonal, meaning that the errors are uncorrelated. When this assumption holds, the optimal solution has a low rank, providing an estimate of the path matrix for the original problem. 

Despite the difference in our estimation formulation and the original one, we believe that our proposed formulations serve two folds. Firstly, unlike previous SEM applications where only a few variables are of interest~\cite{McL:94,BHHB+:00,KZCB+:07,JKCH+:09}, many applications tend to consider a much larger number of variables, such as fMRI studies where the variables are neuronal activities and number up to thousands. Existing approaches of learning causal structures in the exploratory SEM may experience a computational difficulty in terms of the memory storage or convergence. Secondly, our solution for confirmatory SEM is obtained under an assumption of homoskedasticity of residual errors, so if this assumption holds, our and the original solution coincide. Even if it does not hold, and our solution is then not optimal for the original problem, our solution can serve as a starting value for the iterative algorithm used in the original one when a convergence was not obtained.

Related formulations of using SEMs in order to reveal causal structures in multivariate variables chave been reported~\cite{BMG:14,SBG:17,TSG:17,CBG:11}, where the models described a relationship between the output and exogenous variables. These researches applied a sum-of-norm regularization (regarded as an $\ell_1$-type penalty) to promote a sparsity in the path matrix that further inferred a sparse network topology. In~\cite{BMG:14}, the authors proposed a formulation that includes a forgetting factor in the cost objective to track the change in the dynamic network over time, while the contribution of~\cite{SBG:17} was to propose a kernel-based nonlinear SEM model to capture nonlinear features in some applications. One common underlying characteristic in these works is the use of the LS objective as a goodness of fit and an $\ell_1$ norm as a regularization term, so that the resulting formulation was convex and where the variables were the coefficients from exogenous variables to the output. However, it is known that the noise covariance in a LS framework is not a variable that can be estimated and LS estimators are not efficient, in contrast to ML estimators. Our approach exploits the simplicity of a linear SEM but we consider the original SEM formulation that uses Kullback-Leiber (KL) divergence as a goodness of fit while both the path coefficients and the error covariance matrix are estimated simultaneously making the original formulation problem non-convex.  

The two proposed formulations were initially developed in prior work~\cite{PrS:16b, Pru:17} and here we have adapted some details on the algorithms, provided mathematical proofs of important results and included two more experimental real world applications. Section~\ref{sec:sembg} summarizes the mathematical formulation of the path analysis problem as the ML estimation with a quadratic equality constraint. Section~\ref{sec:cvx} describes our convex formulation for confirmatory SEM and shows that the solution can be further used under the condition of having a low rank solution at the optimum. Another convex formulation for exploratory SEM is proposed in Section~\ref{sec:cvxl1}, where an $\ell_1$ regularization is introduced in the cost objective. We show that sparse SEM solutions can be obtained by controlling the regularization parameter. To select the best causality structure, the model selection procedure is explained in Section~\ref{sec:learning_path}.  Then, proximal algorithms for solving the exploratory SEM are explained in Section~\ref{sec:alg}, where we provide some examples of implementation in a large-scale. Numerical experiments in Section~\ref{sec:gendata} demonstrate the important factors to the performance of our method from the simulated data. The results of applying our estimation formulation to learn the causal relations of air pollution data and brain regions are explained in Section~\ref{sec:realdata}. Finally, some mathematical proofs used in our analysis are provided in the Appendix.

\paragraph{Keywords:} structural equation model, convex optimization, regularization, brain connectivity

\paragraph{Notation.} $\symm^n$ denotes the set of symmetric matrices of size $n \times n$ and $\symm^n_{+}$ denotes the set of positive semidefinite matrices of size $n \times n$. For a square matrix $X$, $\Tr(X)$ is the trace of $X$ and $\diag(X)$ is a diagonal matrix containing diagonal entries of $X$. A block symmetric matrix of size $2n \times 2n$ is of the form: $X = \begin{bmatrix} X_1 & X_2^T \\ X_2 & X_4 \end{bmatrix}$, where $X_k$ is of size $n \times n$. The notations $X \succ 0$ and $X \succeq 0$ refer to $X$ being positive definite, and positive semidefinite, respectively. 

\section{Path analysis in SEM}
\label{sec:sembg}
The SEM starts with a set of variables involved in the study, as the measured variables and latent variables. Measured variables are simply the ones that can be directly measured (physical quantities), while latent variables are those that cannot be directly (or exactly) measured, such as intelligence, attitude, etc. Each of these variables can be regarded as either endogenous or exogenous. An endogenous variable gets an influence from others, while an exogenous variable affects the other variables. A general mathematical model in SEM explains a linear relationship from latent variables to measured variables and also includes error terms of each variable~\cite{Bol:89,Mul:09,Hoy:95}.

In practice, we are commonly interested in the application of SEM in situations that only involve \emph{observable} variables. For this reason, we focus on a special class of model in SEM that is described by the multiple linear regression: $Y = c + AY + \epsilon$, where $Y \in \reals^n$ is the measured (or observed) variables, $c \in \reals^n$ is a constant vector representing a baseline and $\epsilon \in \reals^n$ is the model error, assumed to be Gaussian distributed. The \emph{path matrix} $A=[A_{ij}] \in \reals^{n \times n}$ represents a dependence structure among variables in the model. Thus, if $A_{ij} = 0$ then there is no path from $Y_j$ to $Y_i$. In other words, a pattern of nonzero entries in $A$ reveals a causal structure of variables in the model. If this structure is assumed from \emph{prior} knowledge, then the problem of estimating $A$ is called \emph{confirmatory SEM}. 

Let $S$ be a sample covariance matrix of $Y$ and $\Sigma$ be the model-reproduced covariance matrix of $Y$, derived from the regression model, which is then given by  (1);
\begin{equation}
\Sigma = (I-A)^{-1} \Psi (I-A)^{-T},
\label{eq:sigmaAPsi}
\end{equation}
where $\Psi = \Cov(\epsilon)$. The estimation problem in SEM is to find $A$ and $\Psi$ that minimize the KL divergence function $d(S,\Sigma) = \log \det \Sigma + \Tr(S\Sigma^{-1}) - \log \det S - n$, meaning that $\Sigma$ is close to $S$, while maintaining that $\Sigma, A$ and $\Psi$ are constrained by~\eqref{eq:sigmaAPsi}. Moreover, the structure of the path matrix is presumably encoded by a model hypothesis: (i) $A_{ij} = 0$ if there is no link from $Y_j$ to $Y_i$ and (ii) we always have $\diag(A) = 0$, since there is no path from $Y_i$ to itself. To specify the zero structure of $A$, we then define the associated index set $I_A \subseteq \{ 1,2,\ldots,n\} \times \{ 1,2,\ldots,n\}$ with properties that i) $(i,j) \in I_A$ if $A_{ij} = 0$ and ii) $\{(1,1),(2,2),\ldots,(n,n)\} \subseteq I_A$. In short, $I_A$ denotes the index set of \emph{hypothetical zero entries} in $A$. Given an index set $I_A$, we define a \emph{projection operator} $P:\reals^{n\times n} \rightarrow \reals^{n \times n}$, as shown in~\eqref{eq:proj_operator}; 
\begin{equation}
P(X) = \begin{cases} X_{ij} , & (i,j) \in I_A, \\
0, & \mbox{otherwise}, 
\end{cases}
\label{eq:proj_operator}
\end{equation}
and denote $P^c = I-P$. The operators $P^c$ and $P$ are both self-adjoint, \ie, $\Tr(Y^TP(X)) = \Tr(P(Y)^TX)$ and that $P^c(P(X)) = 0$. These two projection operators are then used in the duality of our estimation formulations. With the definition of $P$, the estimation problem corresponding to the confirmatory SEM is shown in~\eqref{eq:sem_original};
\begin{equation}
\begin{array}{ll} 
\minimize & -\log \det \Sigma^{-1}  + \Tr (S \Sigma^{-1}) - \log \det S -n \\
\mbox{subject to } & \Sigma^{-1} = (I-A)^T \Psi^{-1} (I-A), \quad P(A) = 0,
\end{array}
\label{eq:sem_original}
\end{equation}
with variables $A \in \reals^{n \times n}, \Psi \in \symm^{n}_+$ and $\Sigma \in \symm^n_+$. The condition $P(A)=0$ explains the zero constraint on the entries of $A$, and when there is no information on the path matrix, this condition reduces to $\diag(A) = 0$. The problem in~\eqref{eq:sem_original} is one of estimation formulations considered in an SEM context~\cite[\S 4]{Bol:89}. Other cost objectives are also used, such as the ordinary or weighted LS. Estimation formulations in SEM typically consider the ML formulation using the KL divergence objective as in~\eqref{eq:sem_original}, because an ML estimator has the efficient property while it is known that a LS formulation does not always lead to an efficient estimator. 

\paragraph{Special case.} If the constraint $P(A)=0$ reduces to $\diag(A) = 0$ (we allow $A$ to have as many free parameters as possible), then we can make the cost objective zero by solving $S^{-1}= (I-A)^T \Psi^{-1} (I-A)$, where $S$ is given and $A$ and $\Psi$ are free variables. In this case, solutions $A$ of~\eqref{eq:sem_original} are not unique; one can obtain $A$ as dense (non-recursive model) or lower triangular matrix (recursive model). This could be problematic in reading a causality structure from a zero pattern in the estimated $A$. For this reason, it is common to assume some structure in $A$ and a diagonal structure in $\Psi$ (meaning the error terms are uncorrelated). Specifically, the degrees of freedom ($\mathrm{df}$) are defined as in~\eqref{eq:df}, 
\begin{equation}
	\mathrm{df} = \text{ the number of known parameters } - \text{ the number of estimated parameters }.
\label{eq:df}
\end{equation}
Referring to~\eqref{eq:sem_original}, the number of known parameters is the number of entries in the sample covariance matrix and is equal to $n(n-1)/2$ where $n$ is the number of observed variables. The number of free parameters in~\eqref{eq:df} is the total number of entries in $A$ plus the total number of entries in $\Psi$. One can use $\mathrm{df}$ as a guideline for identifying the uniqueness of solution. When the df is negative, the estimator may not be unique. We say that the model is \emph{identifiable} if the df is nonnegative~\cite[p. 35]{RaM:06}. 

\paragraph{Connection with Gaussian graphical models.} Without the parametric covariance constraint in~\eqref{eq:sigmaAPsi}, the cost objective in~\eqref{eq:sem_original} alone can be viewed as a problem of estimating the inverse covariance matrix of Gaussian vectors from a sample covariance $S$. When the objective function is added with an $\ell_1$-norm penalty as $-\log \det X + \Tr (SX) + \sum_{i\neq j} |X_{ij}|$, the problem is known as a \emph{graphical lasso}~\cite[\S 17.3.2]{FHT:09}, which finds many extensions and algorithm developments~\cite{MaH:12,BCP:18}. It is worth nothing that dependence structures among variables decoded from SEMs and the graphical lasso are conceptually different. In graphical lasso, zeros in the inverse covariance matrix $(X)$ of Gaussian vectors explain the conditional independence structure, while no parametric model of variables is assumed. On the other hand, linear SEM assumes a regression model and the interconnection structure of variables are explained from the significant entries of the path matrix $A$, not $\Sigma$. Combining the concept of a conditional independence graph and linear SEM model was proposed in~\cite{LoB:14}, which established a condition of undirected graphs on the inverse covariance $\Sigma^{-1}$, where $\Sigma$ is an implied covariance in~\eqref{eq:sigmaAPsi} from an underlying linear SEM with the assumption that $\Psi$ is diagonal. Their formulation is to estimate $A$ and $\Psi$ so that a conditional independence graph can be recovered. 


\section{Convex formulations for path analysis in SEM}
\label{sec:cvx}
The non-convex property of~\eqref{eq:sem_original} from the quadratic equality constraint is apparent and leads to the chance of obtaining local minima when solving the problem numerically. This section describes the first contribution of our work, where we propose alternative convex formulations and their dual problems for both confirmatory and exploratory SEM. We consider a special case of the path analysis problem where the covariance error is allowed to be diagonal, suggesting that the residual errors of the model were assumed to be uncorrelated. The solution to our formulations is useful when it is low rank at optimum, which was shown to occur under some mild conditions on a problem parameter. The solutions to our formulation and the original problem agree when the covariance of residual error is specified to be a multiple of the identity matrix, which is often the case in SEM applications.

\subsection{Confirmatory SEM}
Define a change of variables as shown in~\eqref{eq:changevar}:
\begin{equation}
X = \begin{bmatrix} X_1 & X_2^T \\ X_2 & X_4 \end{bmatrix},\quad X_1 = \Sigma^{-1}, \quad X_2 = I-A, \quad X_4 = \Psi.
\label{eq:changevar}
\end{equation}
Previous work~\cite{PrS:16b} applied a convex relaxation to the quadratic equality constraint of~\eqref{eq:sem_original} and proposed the \emph{convex confirmatory SEM} formulation
\begin{equation}
\begin{array}{ll}
\minimize & -\log \det X_1 + \Tr(SX_1) \\
\mbox{subject to} & X \succeq 0, \quad  0 \preceq  X_4 \preceq \alpha I , \quad  P(X_2) = I,
\end{array}
\label{eq:sem_primal}
\end{equation}
with variables $X \in \symm^{2n}$, where $S \succ 0 $ and $\alpha > 0$ are given parameters. After relaxing the non-convex constraint, we introduced the inequality constraint $X_4 \preceq \alpha I$ to avoid a trivial solution in~\eqref{eq:sem_primal}, such as $X_4$ can be arbitrarily large and $X_2 = I$. We justify that $\alpha$ can serve as an upper bound on the covariance error of residual in the SEM. Then~\eqref{eq:sem_primal} is a semidefinite programming that can be solved by many existing convex program solvers. It was previously noted~\cite{PrS:16b} that~\eqref{eq:sem_original} and~\eqref{eq:sem_primal} are no longer equivalent, but the convex formulation, that is solved more efficiently, can provide a useful initial solution when solving~\eqref{eq:sem_original} numerically.

Let $Z = \begin{bmatrix} Z_1 & Z_2^T \\ Z_2 & Z_4 \end{bmatrix} \in \symm^{2n}$ be the Lagrange multiplier of the constraint $X \succeq 0$ in~\eqref{eq:sem_primal}. Then the dual problem of the convex confirmatory SEM in~\eqref{eq:sem_primal} is 
\begin{equation}
\begin{array}{ll}
\minimize & -\log \det (S-Z_1) -2 \Tr(Z_2) - \alpha \Tr(Z_4) + n \\
\mbox{subject to} & Z  \succeq 0, \;\; P^c(Z_2) = 0,
\end{array}
\label{eq:sem_dual}
\end{equation}
with variable $Z \in \symm^{2n} $. The constraint $P^c(Z_2) = 0$ explains that the corresponding entries of block $Z_2$ to the zero entries in $X_2$ are free, and the other entries of $Z_2$ are all zero. If the condition $P(X_2) = I$ reduces to $\diag(X_2) = I$ in~\eqref{eq:sem_primal}, then in the dual, it is simplified to that $Z_2$ is diagonal. 

\subsubsection*{Trivial dual solutions}
\label{sec:trivialsol}
The proposed framework in~\eqref{eq:sem_primal} has two problem parameters: $S \succ 0 $ and $\alpha > 0$. (Note that if $S \succeq 0$ the problem could be unbounded below). The important theoretical results that $\alpha$ cannot be arbitrarily large as the solution $X_1$ in~\eqref{eq:sem_primal} becomes an impracticable approximate of $\Sigma^{-1}$ in the original problem of~\eqref{eq:sem_original}~\cite{PrS:16b}. For the completeness of this paper, we state the proposition in~\cite{PrS:16b} and provide the proof here. Let us start with the zero gradient condition of Karush-Kuhn-Tucher (KKT): $X_1= (S-Z)^{-1}$, when $Z=0$ is an optimal dual. Let $\tilde{X}_2 = X_2 -P(X_2) + I$ for any $X_2$. Then the optimal primal solution $X$ must satisfy the feasibility condition in~\eqref{eq:invS_feas},
\begin{equation}
X_1 = S^{-1} ,\quad X_1 \succeq \tilde{X}_2^T X_4^{-1} \tilde{X}_2, \quad 0 \prec X_4 \preceq \alpha I.
\label{eq:invS_feas}
\end{equation}
It suggests that under such acase (trivial dual solution) if $\alpha$ is too large then $X_4$ can be chosen sufficiently large that the choice of $X_1=S^{-1}$ is feasible and optimal, which becomes merely a simple estimate of $\Sigma^{-1}$ and an unfavorable solution. Therefore, the proposition addresses a criterion of the choice of $\alpha$ to avoid such a solution.
\begin{proposition}
\label{prop:alphac}
Let $\alpha_c = n/\Tr(S^{-1})$ (the harmonic mean of the eigenvalues of $S \succ 0$). If $\alpha \leq \alpha_c$ the feasibility problem in~\eqref{eq:invS_feas} has no solution. 
\end{proposition}
See the proof in Appendix~\ref{proof:alphac} that applies Farka's lemma to a semidefinite programming.

To apply the formulation of~\eqref{eq:sem_primal}, when samples of data, $Y_1,Y_2,\ldots,Y_N$ are available, one computes $S$ (the sample covariance of $\{ Y_k\}_{k=1}^N$) and choose $\alpha$, which is as of now, suggested to be less than $\alpha_c$. We can show easily that as one simple choice, the minimum eigenvalue of $S$ denoted as $\lambda_{\mathrm{min}}(S)$ is always less than $\alpha_c$. Suppose $\lambda_1 \leq \lambda_2 \leq \cdots \leq \lambda_n$ are eigenvalues of $S$. It follows that 
$ \frac{1}{\lambda_1} + \frac{1}{\lambda_2 } + \cdots + \frac{1}{\lambda_n} \leq \frac{n}{\lambda_1}$. Since $\Tr(S^{-1}) = \sum_{k=1}^n 1/\lambda_k$, we have
\[
\alpha_c = \frac{n}{\Tr(S^{-1})} = \frac{n}{ \frac{1}{\lambda_1} + \frac{1}{\lambda_2 } + \cdots + \frac{1}{\lambda_n}} \geq {\lambda_1} = \lambda_\mathrm{\min}(S).
\]

\subsection{Sparse SEM with $\ell_1$-norm regularization}
\label{sec:cvxl1}
In exploratory SEM analysis, one aims to discover a zero structure of $A$ from the estimation process that reveals a causal structure of how one variable affects another. One existing approach performs a local search starting from a null model (all path coefficients are zero) and sequentially allows the coefficient corresponding to the largest Lagrangian multiplier to be nonzero~\cite{BHHB+:00}. Another method is to also start from a null model and then add an extra path corresponding to the lowest minimized ML discrepancy function selected among all possible paths. This scheme is referred to as \emph{tree growth} as the model grows by a single entry in $A$ at a time~\cite{CGSH+:11}. The most optimal, but not tractable, approach is to perform a simple brute-force method (or known as \emph{forest growth}) that searches through all possible patterns of zero structures in $A$ with a fixed number of paths and chooses the model corresponding to the lowest minimized ML~\cite{CGSH+:11}. However, the number of all possible models grows exponentially to the number of variables $(n)$, and so it is not feasible as the problem dimension increases.

In this section, we propose a convex formulation for the exploratory SEM problem by applying the widely-used sparse optimization with $\ell_1$-norm. An estimation problem with a $\ell_1$-norm penalty (\emph{lasso}) problem has been well-understood as a convex relaxation of the $\ell_0$-norm penalty that aims to promote a sparsity of the solution~\cite[\S 3]{HTF:09}. This approach has been extended in various directions and applied to find a sparse dependence structure of variables. These include a $\ell_1$-penalty of an affine transformation, termed a generalized lasso~\cite{DSBI+:16}, estimation of the sparse vector autoregressive processes for inferring the Granger causality of time series~\cite{BVN:11}, joint sparse estimation of inverse covariance matrices to find a common conditional independence structure~\cite{THWX+:16, MaM:16}, and estimation of the sparse spectral density matrix to explain the conditional independence structures of time series~\cite{JHG:15, SoV:10}.

Let $h(A) = \sum_{(i,j) \notin I_A} |A_{ij}|$. When $h$ is added to the estimation objective function, it is a form of $\ell_1$-like penalty function (or regularization) as it resembles the $1$-norm of a matrix except, that only those entries not belonging to $I_A$ are penalized. It is well-known that an $\ell_1$-regularized problem returns a sparse solution; many entries of $A$ are zero and the zero pattern can then reveal a causality structure of the variables. Users get to hypothesize about the \emph{known} location of zeros in $A$, which is encoded as the index set $I_A$. Without any prior knowledge about the zero locations in $A$ at all, the constraint $P(A) = \diag(A) =0$ was applied. For $(i,j) \notin I_A$, it is then uncertain if $A_{ij}$ would be zero or not, so the $\ell_1$ norm on these entries is enforced. From the change of variables in~\eqref{eq:changevar}, we have $A_{ij} = -(X_2)_{ij}$ for $(i,j) \notin I_A$. The convex formulation for learning zero patterns in $A$ is then proposed in~\eqref{eq:seml1_primal}; 
\begin{equation}
\begin{array}{ll}
\minimize & -\log \det X_1 + \Tr(SX_1) + 2 \gamma \sum\limits_{(i,j) \notin I_A} |(X_2)_{ij}| \\
\mbox{subject to} & X \succeq 0, \quad 0 \preceq X_4 \preceq \alpha I, \quad  P(X_2) = I,
\end{array}
\label{eq:seml1_primal}
\end{equation}
with variables $X \in \symm^{2n}$. The parameter $\gamma >0$ is the \emph{regularization parameter} that controls the sparsity of $X_2$ (and hence, in $A$). Let $Z = \begin{bmatrix} Z_1 & Z_2^T \\ Z_2 & Z_4 \end{bmatrix} \in \symm^{2n}$ be the Lagrange multiplier of $X \succeq 0$ in~\eqref{eq:seml1_primal}. The dual of~\eqref{eq:seml1_primal} is
\begin{equation}
\begin{array}{ll} \mbox{maximize} & \log \det (S-Z_1) -2 \Tr(Z_2) - \alpha \Tr(Z_4) +n , \\
\mbox{subject to} & Z \succeq 0, \quad \Vert P^c(Z_2) \Vert_\infty \leq \gamma ,
\end{array}
\label{eq:seml1_dual}
\end{equation}
with variable $Z \in \symm^{2n}$. The difference of the dual in~\eqref{eq:seml1_dual} from that in~\eqref{eq:sem_dual} is that $Z_{ij}$ are not forced to be zero, but allowed to be less than $\gamma$ for $(i,j) \notin I_A$.

\subsubsection*{Choice of penalty parameter}
We see that the sparsity of $X_2$, and hence, of the optimal $A$ can be controlled via $\gamma$. For example, the larger $\gamma$, the sparser $A$ is. In Appendix~\ref{app:gammamax}, we established that there exists a critical value of $\gamma$, denoted by $\gamma_\mathrm{\max}$, in the sense that if~\eqref{eq:gammamax} holds;
\begin{equation}
\gamma \geq \gamma_\mathrm{\max}: = (1/\alpha) \Vert P^c (\alpha I - S) \Vert_\infty
\label{eq:gammamax}
\end{equation}
then the optimal $A$ in~\eqref{eq:seml1_primal} is the zero matrix. Moreover, the value of $\gamma_\mathrm{\max}$ can be calculated in advance and depends only on $\alpha$ and $S$ (problem parameters). This means it is unnecessary to vary $\gamma$ arbitrarily in the problem, and so $\gamma_\mathrm{\max}$ can be used as an upper bound of the range of $\gamma$ used for varying the sparsity patterns of $A$. 

\subsubsection*{Related sparse SEM formulations}
The goal of learning sparse coefficients in SEM is commonly found in recent literature. It is worth noting several related model structures that include exogenous $(X)$, mediate $(M)$ or latent variables $(\eta)$ and describe formulations that encourage zeros in the coefficients of these variables.  The first related formulation is the model containing both dependent and latent variables~\cite{ZTT:16, HCW:17} as $Y = \Lambda \eta + \epsilon_y$ and $\eta = \Pi \eta + \epsilon_\eta$. They aim to promote a sparsity in $\Pi$ (coefficients of latent variables) where  \cite{ZTT:16} considered smoothly clipped absolute deviation (SCAD) formulation, while~\cite{HCW:17} applied various sparse penalties such as SCAD, minimax concave penalty and $\ell_1$-norm. The log-likelihood function certainly contains $\eta$, which is unobservable, so the ECM algorithm was employed to find a local maxima. Examples in \cite{JGM16, ZTT:16, HCW:17} illustrate the case but with a small number of variables (up to $20$ in $Y$ and less than $10$ in $\eta$). This framework is very useful when the latents are considered in the model, while we believe that our path analysis problem (that contains only observable variables) is not precisely a special case of the model in~\cite{ZTT:16, HCW:17} as they do not include a path matrix from $Y$ to $Y$. Therefore, path analysis should be specifically solved where more efficient algorithms can be applied as a result of problem simplification and relaxed convexity. The second relevant work is the inclusion of mediate variables~\cite{SJBG:17} as $Y = cX + \diag(b)M + \epsilon_y$ and $M = a^T X + \epsilon_m$, where $Y$ and $X$ are scalars and $M$ is a multivariate mediator. The aim was to obtain a sparsity in $b$ to regard only significant mediators and so they formulated the models into a path analysis so that \texttt{regsem} using the $\ell_1$-norm penalty in~\cite{JGM16} could be applied. Examples are illustrated with the model having only five mediators. It is noted that their model regards $Y$ and $X$ as single variables. Moreover, when formulated as the path analysis problem: $Z =A Z + \epsilon$ with $Z = (Y,X,M)$, it could also be solved using our framework with a minor modification on the definition of $I_A$ as a block of $A$ is restricted to $I$. The third model, by~\cite{CBG:11}, is described as $Y = BY + FX + \epsilon$ with the goal of estimating a sparse $B$. It is then obvious that our formulation is essentially their special case. However, the LS framework is chosen in~\cite{CBG:11} where noise covariance was not part of the variables and so their formulation is simplified and convex. Similarly, if both $X$ and $Y$ are observable, this model can be formulated as a path analysis with $Z = (Y,X)$ and can be solved using our method.

\subsection{Low rank solution}
\label{sec:lowranksol}
Solutions of the convex confirmatory SEM and the sparse SEM are useful for the original SEM problem if $X_1 = X_2^T X_4^{-1} X_2$ at optimum, as we can then use $X_1$ as an estimate of $\Sigma^{-1}$. This occurs if and only if $\Rank(X)= n$ at optimum (note that $X \in \symm^{2n}$). Therefore, we aimed to find a relation between $\alpha$ and low rank optimal solutions of~\eqref{eq:sem_primal} and~\eqref{eq:seml1_primal} from the optimality conditions. The result in section~\ref{sec:trivialsol} suggested that if $\alpha$ is too large, then it is possible that $\Rank(X) > n$, which is to be avoided. 

Previous analysis~\cite{PrS:16b} showed that when an optimal $X$ has rank $n$, then 
\[
X_1 = X_2^{T}X_4^{-1}X_2 ,\quad \Rank(Z) = n,\quad Z_4 \succ 0, \quad X_4 = \alpha I.
\]

\begin{figure}[ht]
\centering
\includegraphics[width=0.95\linewidth]{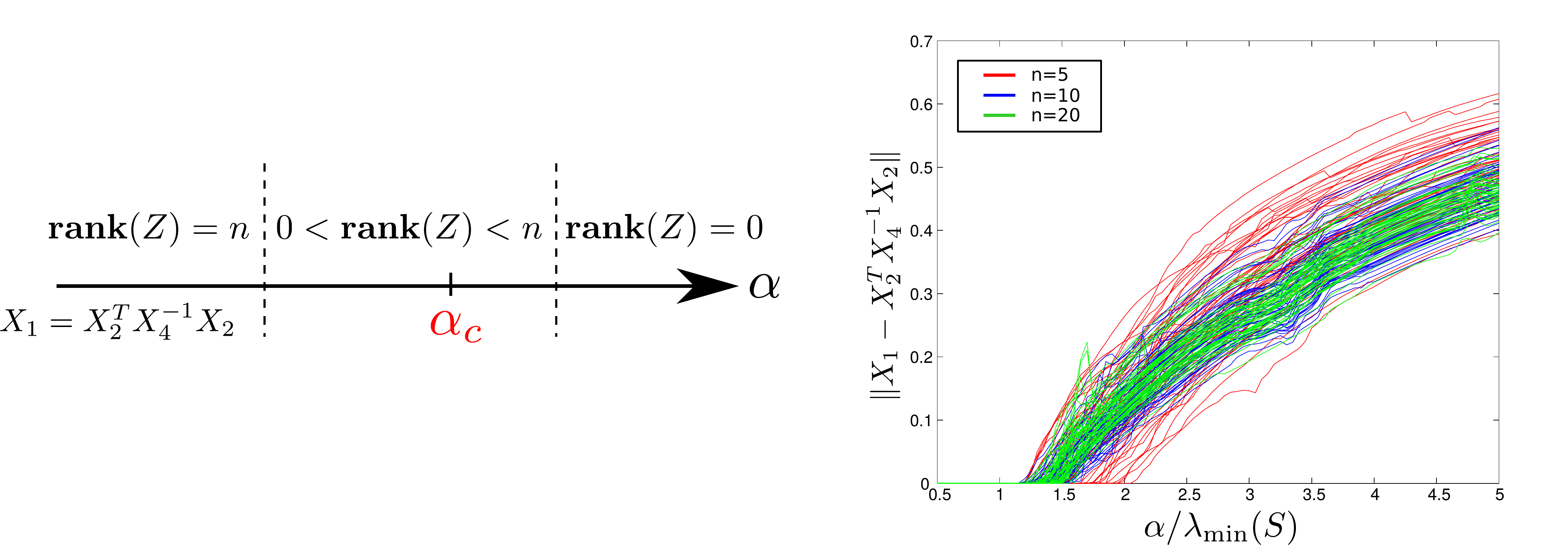}
	\caption{The effect of $\alpha$ on the solutions of the convex confirmatory SEM. \emph{Left.} $\Rank(Z)$ as $\alpha$ varies. \emph{Right.} Error between $X_1$ and the low rank approximation as $\alpha$ varies. Lines with the same color correspond to the result from using the same $n$ and shown as results of many trials.  \text copyright [2016] IEEE. Reprinted, with permission, from the Proceedings of 2016 55th Annual Conference of the Society of Instrument and Control Engineers of Japan (SICE).}
\label{fig:lowrank}
\end{figure}

In Section~\ref{sec:trivialsol} it was shown that the minimum eigenvalue of $S$ always lies on the left of the harmonic mean of the eigenvalues of $S$. If we choose $\alpha = \lambda_\mathrm{\min}(S)$, then it is often the case that $\Rank(Z) = n$. This advises us to include a constraint $X_4 \preceq \lambda_\mathrm{\min}(S)I \preceq S$ into the estimation problem. Since $X_4 = \Psi$, we can justify that the covariance error is controlled to be less than the covariance of the variables. Figure~\ref{fig:lowrank} illustrates that the error between $X_1$ and $X_2^T X_4^{-1} X_2$ increases as $\alpha$ increases and is zero when $\alpha$ is sufficiently small relatively to the minimum eigenvalue of $S$. It is also noted that when $\alpha$ is too large, then $\Rank(X) > n$ and there could be many optimal solutions $X_4$ that is strictly less than $\alpha I$. In this case, the solution $X_1$ is not unique.

\subsection{Exploratory SEM}
\label{sec:learning_path}
As seen in Section~\ref{sec:cvxl1}, controlling the regularization parameter in the sparse SEM problem can provide path matrix solutions with various sparsity patterns. If $\gamma$ is large then the path matrix $A$ contains many zeros, resulting in a simple interpretation of the estimated causal structure but the goodness of fit becomes smaller. Therefore, choosing an appropriate value of $\gamma$ is a trade-off between choosing the solution to explain a causal structure in a simple way and to best describe the data to a certain level. Two common approaches of selecting the regularization parameter implemented in many works on sparse learning are the cross-validation and model selection criterions~\cite{FHT:09}, and examples of applying these approaches in a graphical lasso problem have been reported~\cite{YuL:07,GLMZ:11}. In a $K$-fold cross validation, we vary $\gamma$ over a specified range and train the model on the training data fold, then the $\gamma$ value that yielded the smallest test prediction error~\cite{FHT:09} is chosen. The cross validation approach can be computationally intensive since the models are estimated from the data from $K$-folds and for several values of $\gamma$. 

We adopted a model selection approach in this paper, which is a widely and commonly used approach in many regularized SEM problems~\cite{HCW:17, ZTT:16, SJBG:17}, where the value of $\gamma$ corresponding to the model with the lowest score in the model selection criterion~\cite[\S 7]{FHT:09} was chosen. These included the Akaike Information Criterion (AIC), the corrected AIC (AICc), Bayesian Information Criterion (BIC)~\cite[\S 7]{FHT:09}, Kullback Information Criterion (KIC)~\cite{Cav:99} and the corrected KIC (KICc)~\cite{SeB:04}, where we calculated the log-likelihood function of samples $Y_1,Y_2,\ldots,Y_N$ as $\mathcal{L}  = \frac{N}{2} \left (- \log\det{\hat{\Sigma}} - \Tr(S\hat{\Sigma}^{-1}) \right ) $. These scores consisted of the two terms of the negative log-likelihood, representing the goodness of fit, and the model complexity. Using a large value of $\gamma$ yielded a sparse model, so the model complexity is low but the fitting is worse. Therefore, a fair comparison among models with different sparsity patterns can be performed using these model selection criterions. These scores have different asymptotic behaviors. The BIC is known to prefer a simpler model, since the penalty on the model complexity is higher relatively to other criterions, while the AIC is the first model selection criterion, which is widely accepted and is known to perform poorly when the number of sample sizes is small relative to the number of effective parameters. The AICc and KICc scores were then developed to reduce the bias and improve the model selection in a small-sample setting.  

\begin{figure}[!ht]
\centering
\centerline{\includegraphics[width=0.5\columnwidth]{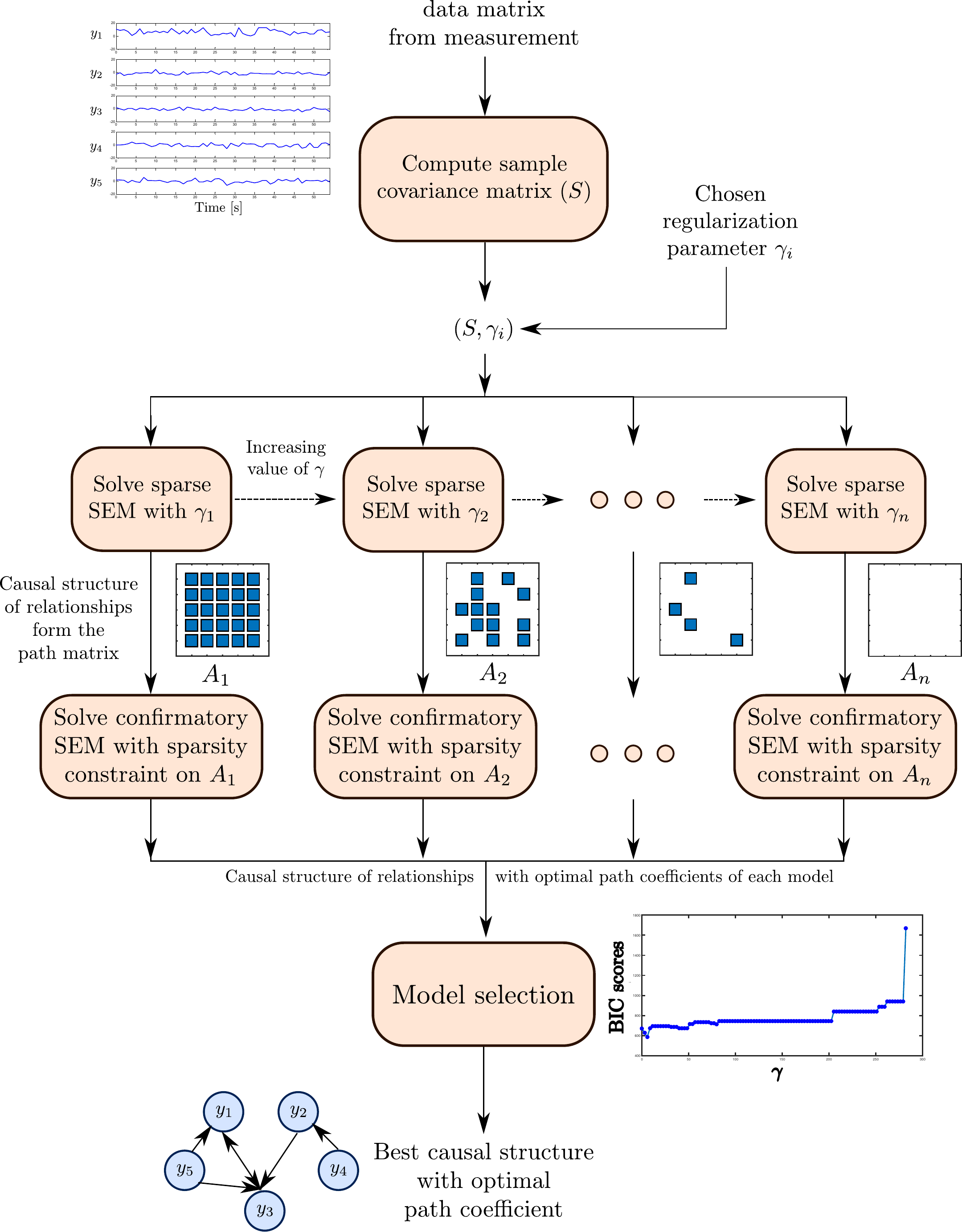}}
\caption{Procedures of learning causal structure of path matrices in SEM. The best optimal causal structure of path matrices is chosen from a model that has the minimum model selection score.}
\label{fig:model_sel1}
\end{figure}

To learn the best causal structure of path matrices, we can choose a range of $\gamma$ values from $0$ to $\gamma_{\max}$ and then solved the sparse SEM in~\eqref{eq:seml1_primal} for each of those values, resulting in estimated path matrices with different sparsity patterns ranging from the densest to the sparsest. Each of the estimated sparsity patterns was then used as a sparsity constraint on $A$ in the convex confirmatory SEM in~\eqref{eq:sem_primal} and was solved for the optimal path matrix $A$ to yield a candidate model. Estimation of the unrestricted model with the selected entries of $A$ by sparse SEM also helps to reduce the bias of the estimator~\cite[\S 3.8.5]{FHT:09}. This process was repeated for all the values of $\gamma$ to obtain a set of candidate models, each of which was labeled by a model selection criterion score and the best model was the one with the minimum score. In short, the sparse SEM was used to select a finite number of sparsity patterns in $A$ and the convex confirmatory SEM was used to provide the best estimate of the path matrix corresponding to the sparsity pattern selected from the model selection criterion score. This procedure is illustrated in Figure~\ref{fig:model_sel1}.

\section{Proximal algorithms}
\label{sec:alg}
In this section, efficient numerical methods were applied to solve the two estimation formulations; the convex confirmatory SEM and the sparse SEM. We examine two proximal algorithms, namely, parallel proximal algorithm (PPXA) and alternating direction method of multipliers (ADMM), that requires a splitting technique on the cost objective of convex problems and introduces some auxiliary variables. Mostly, the key success of proximal algorithms is to add indicator functions corresponding to nonlinear constraints, where the update steps that require a computing of proximal operators can often turn into an almost analytical form. The ADMM algorithm is also examined in related problems, such as an estimation of dynamic SEM~\cite{BMG:14}, kernel-based SEM~\cite{SBG:17} or a SEM with an elastic net~\cite{TSG:17}. In contrast to our ML formulation, all of these works have used a LS objective, so no positive definite constraint of $\Psi$ (error covariance) was needed. However, the algorithms applied to our problem are more involved due to the cone constraint. The main proximal operators involved in our problem are the projection operator of $v$ on a set $\mathcal{C}$. We follow the standard format and implementation guidelines in~\cite{CoP:11,BPCP+:10}. Moreover, the notation of the proximal operator: $\prox_{f}(y) = \argmin_{x} f(x) + (1/2) \Vert y - x \Vert_2^2$ was applied here~\cite{PaB:13}. Applying the algorithms to both the confirmatory and sparse SEMs is quite similar, so only the algorithm detail for the sparse SEM is presented.

\subsection{Algorithms for confirmatory SEM}
\label{sec:alg_confirm}
The confirmatory SEM~\eqref{eq:sem_primal} has three constraints. We split the objective function into two terms by defining functions $f : \symm^{2n} \rightarrow \reals$, $g : \symm^{2n} \rightarrow \reals$ that cooperate the constraints in~\eqref{eq:sem_primal} as indicator functions as
\begin{equation*}
f(Z) = -\log\det(Z_1) + \Tr(SZ_1) + I\{ P(Z_2) = I \} + I \{ 0 \preceq Z_4 \preceq \alpha I\}, \qquad 
g(Z) = I \{ Z \succeq 0 \}.
\end{equation*}
We denote $I \{ Z \in C \}$ as an indicator function with a function value of $0$ if $Z \in C$ and $\infty$ otherwise.  
The ADMM format of the problem in~\eqref{eq:sem_primal} is to $\minimize  f(X) + g(Z)$ subject to $X - Z = 0$ with variables $X , Z \in \symm^{2n}$. For PPXA, we extend the variable into $X= (X_1,X_2)$ where $X_1 = X_2$. The iteration update of both algorithms required deriving the proximal operators of $f$ and $g$, which could be obtained in a similar way to the updates in the sparse SEM problem. Hence, we omit the details here and present them in the next section (section~\ref{sec:alg_sparse_sem}) instead.

\subsection{Algorithms for sparse SEM}
\label{sec:alg_sparse_sem}
In the sparse SEM, we extend the variable into a larger product space and split the cost objective function into the three terms, where each of them corresponds to a simple proximal operator. Let us define functions $f_i:\symm^{2n} \rightarrow \reals$ for $i=1,2$ and 3, as given by:
\begin{gather*}
f_1(Z) = -\log\det(Z_1) + \Tr(SZ_1) + I \{ 0 \preceq Z_4 \preceq \alpha I\} ,\;
 f_2(Z) =  2\gamma \sum_{(i,j) \notin I_A}\vert (Z_2)_{ij}| +  I\{ P(Z_2) = I \} ,\; f_3(Z) = I \{ Z \succeq 0 \}.
\end{gather*}
Then it is obvious that the problem of~\eqref{eq:seml1_primal} is equivalent to the minimization of $f_1+f_2+f_3$. If we extend the variable into a product space then the problem is also equivalent to the form of a \emph{global consensus problem}, shown in~\eqref{eq:consensus};
\begin{equation}
\minimize_{(X_1,X_2,X_3) \in \symm^{2n} \times \symm^{2n} \times \symm^{2n}} \;\; f_1(X_1) + f_2(X_2) + f_3(X_3) \quad \text{subject to} \;\; X_1 = X_2 = X_3,
\label{eq:consensus}
\end{equation}
which can be solved by many existing algorithms, such as the parallel proximal algorithm (PPXA)~\cite[\S 10]{CoP:11}, and ADMM~\cite{BPCP+:10}. The update rules of these two algorithms are described in the Appendix~\ref{sec:alg_description} and require proximal operators of $f_1,f_2$ and $f_3$ as follows.

\paragraph{Proximal operator of $f_1$.} Finding $\prox_{f_1/\rho}(Y)$ involves solving the problem:
\[
\minimize_{Z \in \symm^{2n}} \;\;-\log\det(Z_1)+\Tr(SZ_1) + (\rho/2) \Vert Y - Z \Vert_F^2 \quad \text{subject to} \;\;0 \preceq Z_4 \preceq \alpha I. 
\]
The blocks $Z_1,Z_2$ and $Z_4$ in $Z$ can be optimized independently because the Frobenious norm term is separable. Firstly, $Z_2 = Y_2$ and $Z_4$ is obtained by projecting the block $Y_4$ into the set $\mathcal{C} = \{ Z \;|\;0 \preceq Z_4 \preceq \alpha I \}$. To find $Z_1$, we derive the zero gradient condition of the objective function with respect to $Z_1$ as shown in~\eqref{eq:minX_grad_explor};
\begin{equation}
\rho Z_1 - Z_1^{-1} = \rho Y_1 - S,
\label{eq:minX_grad_explor}
\end{equation}
with an implicit constraint: $Z_1 \succ 0$. To solve for $Z_1 \succ 0 $ that satisfies~\eqref{eq:minX_grad_explor}, we can apply the technique from~\cite[\S 6.5]{BPCP+:10}. We compute the eigenvalue decomposition: $\rho Y_1 - S  = Q\Lambda Q^T$, where $Q$ is the eigenvector matrix and $\Lambda$ is diagonal and contains the eigenvalues. Then, we define $\tilde{Z}_1$ as a diagonal matrix with decomposition $Z_1 = Q\tilde{Z}_1 Q^T$. Substituting these decompositions into~\eqref{eq:minX_grad_explor} gives $\rho z^2- \lambda z - 1 = 0$ where $z$ is the diagonal entries of $\tilde{Z}_1$. We can solve the quadratic equation to obtain the positive root as $z = (\lambda + \sqrt{\lambda^2 + 4 \rho})/2\rho$ and then $Z_1$ is guaranteed to be positive definite. 

\paragraph{Proximal operator of $f_2$.} In this step, to find $\prox_{f_2/rho}(Y)$, we are required to solve the problem of the form:
\[
\minimize_Z \;\; 2\gamma\sum_{(i,j) \notin I_A}\vert Z_{ij} \vert  + (\rho/2)\Vert Y - Z\Vert_F^2,\quad \text{subject to}  \;\;P(Z_2) = I.
\]
The problem is separable again in blocks of $Z$ and the optimized $Z_1$ and $Z_4$ are simply $Y_1$ and $Y_4$, respectively. It is left to find the block $Z_2$ where the entries of $(i,j) \notin I_A$ can be solved from a simple problem: $\minimize_z \;\; 2\gamma|z| + \rho|z-y|^2$, typically known as finding a proximal operator of $f(x) = |x|$. The solution of minimizing the above problem can be performed by element-wise soft thresholding, defined by $(Z_2)_{ij} = S_{\gamma/\rho}((Y_2)_{ij})$, where the \emph{soft thresholding} operator~\cite{BPCP+:10} is defined by:
\[
S_k(a) =  a-k, \quad \text{for}\;\; a > k, \quad S_k(a) = 0 , \quad \text{for}\;\; \vert a \vert \leq k, \quad S_k(a) =  a+k, \quad \text{for}\;\;  a<-k,
\]
or equivalently $S_k(a) = \max(|a|-k)\cdot \sign(a)$. Lastly, for $(i,j) \in I_A$ we must have $P(Z_2) = I$.

\paragraph{Proximal operator of $f_3$.} This step solves the problem: $\minimize_{Z \succeq 0} \;\; \Vert Y - Z \Vert_F^2$. The solution is, therefore, the projection of $Y$ onto the positive definite cone.

\subsection{Choice of algorithm parameters}
\label{sec:algparam}
The choice of ADMM parameter, $\rho$, affects the convergence speed and so an optimal selection is still an open question. Recent progress includes an adaptive formula of $\rho$ varying upon the primal and dual residuals~\cite{BPCP+:10} or the rule motivated by the Barzilai-Borwein spectral method~\cite{ZFG:17}. However, this paper did not follow this direction but rather an experimentally derived practical choice of $\rho$ that led to a convergence by a reasonable number of iterations was used. When considering~\eqref{eq:seml1_primal}, we found that good choices of $\rho$ varied with the problem scale, as determined by $(S,\alpha)$. When the minimum eigenvalue of $S$ was very small~$(\leq 10^{-3})$, the ADMM seemed to converge slowly for $\rho = 10,100$. If the problem of~\eqref{eq:seml1_primal} was scaled properly by $\beta$, the choice of $\rho = \beta$ worked well in practice. 

The PPXA parameters described in~\cite[\S 10]{CoP:11} are $\gamma > 0$ and $\lambda \in (0,2)$. We found in the experiment that if the problem is scaled to have $\lambda_{\min}(S) = 1$, the choice of $\gamma = 0.1$ and $\lambda = 1.8$ yield comparatively smaller number of iterations in our problem.  Hence, from the two observations of implementing the two algorithms, we needed a result where solutions from the scaled and original problems could be interchangeably obtained.

\begin{proposition}
\label{prop:scaled_sem}
Let $(X,Z)$ be the primal and dual optimal solutions of~\eqref{eq:seml1_primal} and~\eqref{eq:seml1_dual}, respectively, using the problem parameter $(S,\alpha)$. Moreover, let $(\tilde{X},\tilde{Z})$ be the primal and dual optimal solutions of~\eqref{eq:seml1_primal} and~\eqref{eq:seml1_dual}, respectively, using the problem parameter $(\tilde{S},\tilde{\alpha})$, where $\tilde{S} = \beta S$ and $\tilde{\alpha} = \beta \alpha$. Then, $(X,Z)$ and $(\tilde{X},\tilde{Z})$ are related by~\eqref{eq:scaled_XZ},
\begin{equation}
\tilde{X} = \begin{bmatrix} X_1/\beta & X_2^T \\ X_2 & \beta X_4 \end{bmatrix},\quad \tilde{Z} = \begin{bmatrix} \beta Z_1 & Z_2^T \\ Z_2 & Z_4/\beta \end{bmatrix}.
\label{eq:scaled_XZ}
\end{equation}
In other words, we can obtain the optimal solution of the scaled sparse SEM instantly from the solution of the unscaled problem, and vice versa.
\end{proposition}
The proof of this result is described in Appendix~\ref{sec:scaled_sem}. To apply Proposition~\eqref{prop:scaled_sem}, we scale the sparse SEM by $\beta = 1/\lambda_{\min}(S)$, so that the minimum eigenvalue of $\tilde{S} = \beta S$ is one and $\tilde{\alpha} = 1$. As a result, we solve the scaled sparse SEM using $(\tilde{S},\tilde{\alpha})$ with a heuristic choice of parameters described in section~\ref{sec:algparam} and retrieve the solution to the original problem using~\eqref{eq:scaled_XZ}.

\subsection{Algorithm performance}
\label{sec:algo_perf}
We show the performance of PPXA algorithm in solving the sparse SEM in~\eqref{eq:seml1_primal} with $n = 100,200,\dots,500$ using 50 samples of $S$ for each $n$. We set $20\%$ of entries in $A$ to satisfy the constraint $P(A) = 0$ in random locations. The algorithm performances generally depended on the optimization problem parameters. In this experiment, we generated $S$ as a covariance matrix of the $n$-dimensional random vector $Y$ and varied the eigenvalues of $S$ into small $\lambda_{\min}(S) \in  (0.08,0.09)$ or large $\lambda_{\min}(S) \in (0.47,0.48) $ using $2n$ and $10n$ samples of $Y$, respectively. A slower convergence was expected when $\lambda_{\min}(S)$ was smaller. The second factor in a convergence is the regularization parameter, so we varied $\gamma$ as $0.05 \gamma_{\max}$ (dense solution) or $0.8 \gamma_{\max}$ (sparse solution). The algorithm performances were evaluated in terms of the number of iterations required to satisfy the stopping criterion and by the CPU times which are averaged over $50$ runs. The program stops when the relative change of the cost objective and the relative change of the solution are less than $10^{-5}$ (low accuracy). The  specifications of the computer used in this experiment are: CPU: Intel Core I7-8559U (4.5 GHz), RAM: 16GB DDR4 BUS2400, HD: SSD 500GB, OS: Ubuntu 18.04LTS. The codes were programmed in MATLAB and parallel computing toolbox was used in this simulation. 

\begin{figure}[ht]
\centering
\includegraphics[width=0.8\linewidth]{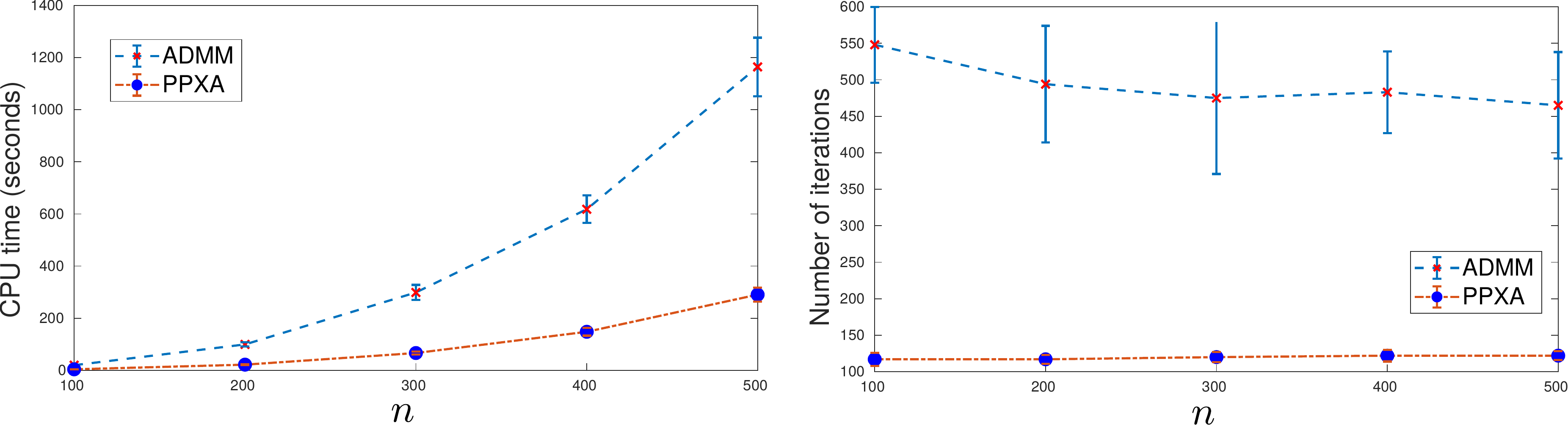}
\caption{Algorithm performances as $n$ varies under the setting of small $\lambda_{\min}(S)$ and $\gamma = 0.05\gamma_{\max}$.}
\label{fig:alg_admm_ppxa}
\end{figure}

A comparison of convergence between ADMM and PPXA are shown in Figure~\ref{fig:alg_admm_ppxa} where PPXA appears to take less number of iterations (and hence less CPU time) than ADMM to achieve the same accuracy.
Hence, only results of PPXA algorithm are reported in Tables~\ref{tab:ppxa_time} and~\ref{tab:ppxa_iter}. Table~\ref{tab:ppxa_iter} shows that when $S$ was almost degenerated, a higher number of iterations was required to achieve a desired accuracy. In addition, the algorithm converges slightly slower when using a large $\gamma$ to obtain sparse solutions but the effect of varying $\gamma$ was not apparent when $\lambda_{\min}(S)$ was large. Solving the confirmatory and sparse SEMs with dimension $n$ involves the total number of variables in $X$. The main computational cost of the two proximal algorithms for both problems depend on the eigenvalue decomposition in a symmetric matrix of size $2n$, which is ${\mathcal O}((2n)^3)$. A trial problem with $n = 500$ and a given pattern in $A$, resulting in a total of $250,000$ variables, required approximately $4-9$ minutes of computational time in low accuracy mode. Such a large-scale setting may not be feasible when applying an iterative method based on the use of a Hessian matrix.  

\begin{table}[!ht] 
     \centering 
     \caption{Averaged CPU times and standard deviations (in parentheses) of PPXA algorithm.} 
     \begin{tabular}{|c| c|c|c|c|} 
     \hline 
     \multirow{3}{*}{$n$}  & \multicolumn{4}{|c|}{CPU time (seconds)}\\ 
     \cline{2-5} 
     & \multicolumn{2}{|c|}{$\gamma = 0.05\gamma_{\mathrm{max}}$} & \multicolumn{2}{|c|}{$\gamma = 0.8\gamma_{\mathrm{max}}$}\\ 
     \cline{2-5} 
     & small $\lambda_{\min}(S)$ & large $\lambda_{\min}(S)$ & small $\lambda_{\min}(S)$  &large $\lambda_{\min}(S)$ \\ 
     \hline 
     100  & 3.92 (0.33)  & 3.14 (0.19) & 6.19 (0.62) &  3.42 (0.28) \\ 
     200  & 21.83 (1.81)  & 17.73 (1.24) & 35.78 (2.89) &  20.06 (1.64) \\ 
     300  & 66.41 (5.52)  & 52.74 (3.95) & 109.24 (8.65) &  60.63 (5.46) \\ 
     400  & 147.94 (15.05)  & 115.42 (8.96) & 245.30 (22.07) &  134.58 (11.14) \\ 
     500  & 290.48 (26.74)  & 222.03 (16.18) & 495.80 (41.11) &  267.63 (25.54) \\ 
     \hline 
     \end{tabular} 
     \label{tab:ppxa_time} 
     \end{table} 

\begin{table}[!ht] 
     \centering 
     \caption{Averaged number of iterations and standard deviations (in parentheses) of PPXA algorithm.} 
     \begin{tabular}{|c| c|c|c|c|} 
     \hline 
     \multirow{3}{*}{$n$}  & \multicolumn{4}{|c|}{Number of iterations}\\ 
     \cline{2-5} 
     & \multicolumn{2}{|c|}{$\gamma = 0.05\gamma_{\mathrm{max}}$} & \multicolumn{2}{|c|}{$\gamma = 0.8\gamma_{\mathrm{max}}$}\\ 
     \cline{2-5} 
     & small $\lambda_{\min}(S)$ & large $\lambda_{\min}(S)$ & small $\lambda_{\min}(S)$  &large $\lambda_{\min}(S)$ \\ 
     \hline 
      100 & 117 (9)  & 93 (2) & 215 (19) &  112 (5) \\ 
      200 & 117 (6)  & 92 (2) & 221 (12) &  116 (3) \\ 
      300 & 120 (7)  & 92 (2) & 225 (10) &  117 (2) \\ 
      400 & 122 (8)  & 91 (3) & 227 (11) &  118 (2) \\ 
      500 & 122 (6)  & 90 (3) & 226 (8) &  118 (2) \\ 
     \hline
     \end{tabular} 
     \label{tab:ppxa_iter} 
     \end{table} 
 
We note that it is not our goal to make an intensive comparison of available algorithms in this study. Both algorithms are available for publicly tested at \texttt{https://github.com/jitkomut/cvxsem}.

\section{Results (generated data)}
\label{sec:gendata}
This section describes the performance of the exploratory SEM presented in Section~\ref{sec:learning_path} with the scheme in Figure~\ref{fig:model_sel1}. The goal was to examine how well the nonzero and zero entries in the path matrix could be estimated from the data generated from the true model encoded with a pattern of the true path matrix. Throughout this section, we set $n=10$ and choose $A_{\mathrm{true}}$ having random sparsity patterns and generated measurements from $ Y = (I - A_{\mathrm{true}})^{-1}e$, where $e$ is normally distributed with a variance of $0.1$. The matrix $S$ was then computed as the sample covariance of $Y$. The options of the data generating process were the number of samples and the density of non-zero elements in $A_\mathrm{true}$. In the estimation process, as required by the constraint $P(A) = 0$, we made assumptions about zero locations in $A$, varied by patterns and as the percentage of all zeros. To examine the performance, we applied the typical measures of true positives (TP), true negative (TN), false positives (FP), false negatives (FN), TP rate and FP rate, where positives are the non-zero entries in $A$ and the negatives are zeros in $A$, through a receiver operating characteristic (ROC) curve~\cite{Alp:14}[\S 19.7].
\subsection{Performance of the sparse SEM}
Four main aspects that could influence the estimation results were explored. These factors were the sparsity density of the true model ($A_{\mathrm{true}}$), the number of samples $(N)$, the number of known zero locations used in the estimation and the choice of model selection scores. The experiments were then designed to investigate the effects of these factors as explained below.

\begin{figure}
\centering
    \begin{subfigure}[b]{0.9\linewidth}            
            \includegraphics[width=\textwidth]{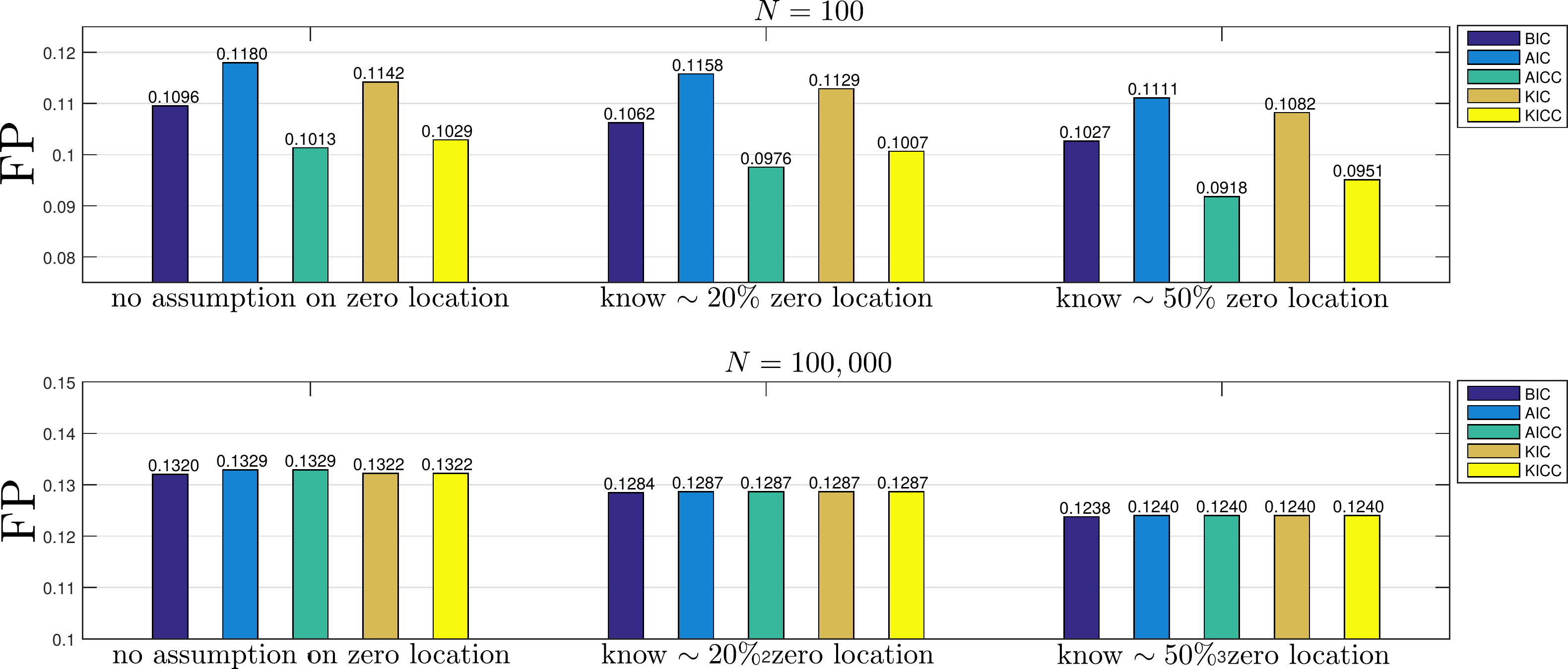}
            \caption{\textbf{False Positive (FP) error.} When $A_{\mathrm{true}}$ is dense, the FP from all model selection criterions tended to decrease with increasing knowledge about the zero location in $A_{\mathrm{true}}$, but increased when $N$ was high as all model selection criterions tended to select the denser $\hat{A}$. In the case of a small $N$, AICc provided the minimum FP error.}
            \label{fig:FP_Atrue_dense}
    \end{subfigure}%

    \begin{subfigure}[b]{0.9\linewidth}
            \includegraphics[width=\textwidth]{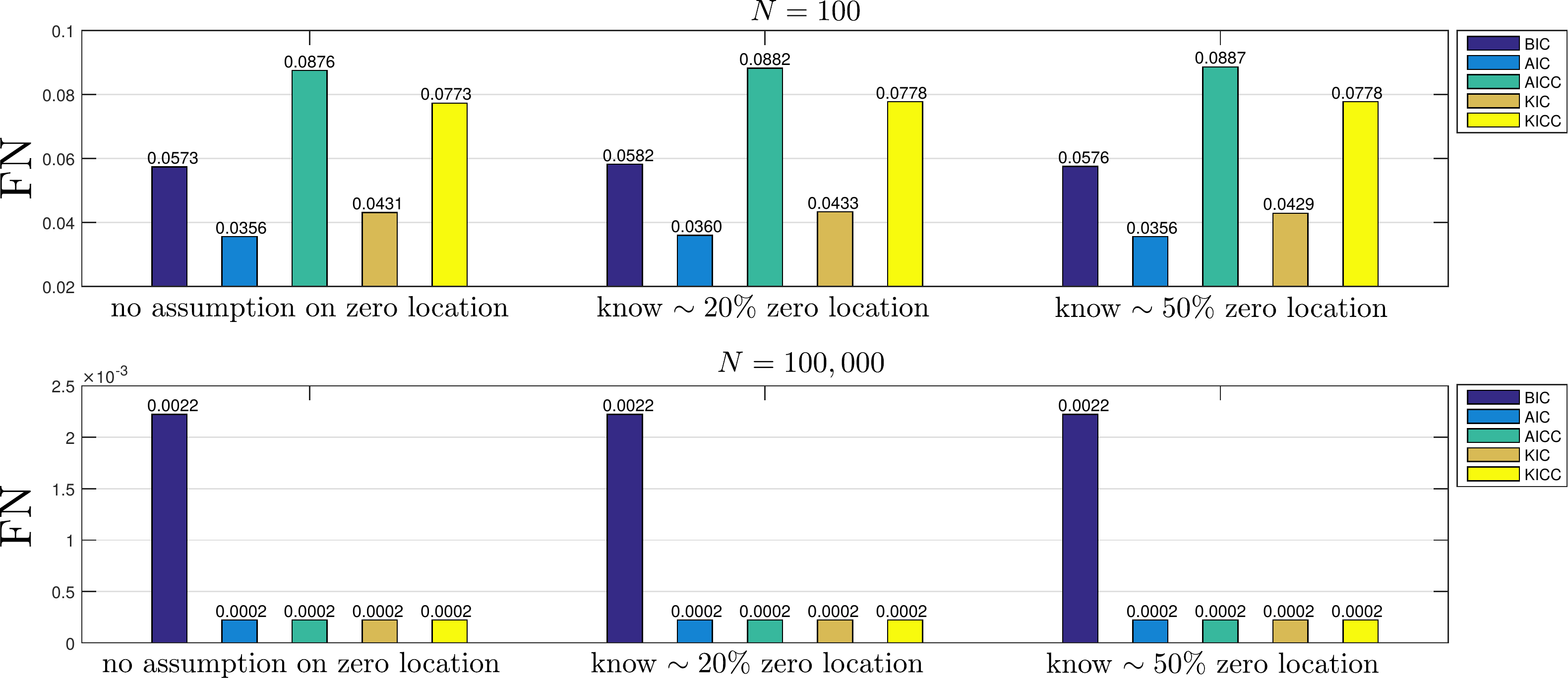}
            \caption{\textbf{False Negative (FN) error.} When increasing the knowledge about the zero location in $A_{\mathrm{true}}$ the estimation process barely affected the change of FN, but it was improved when $N$ grows.}
            \label{fig:FN_Atrue_dense}
    \end{subfigure}

    \begin{subfigure}[b]{0.9\linewidth}
            \includegraphics[width=\textwidth]{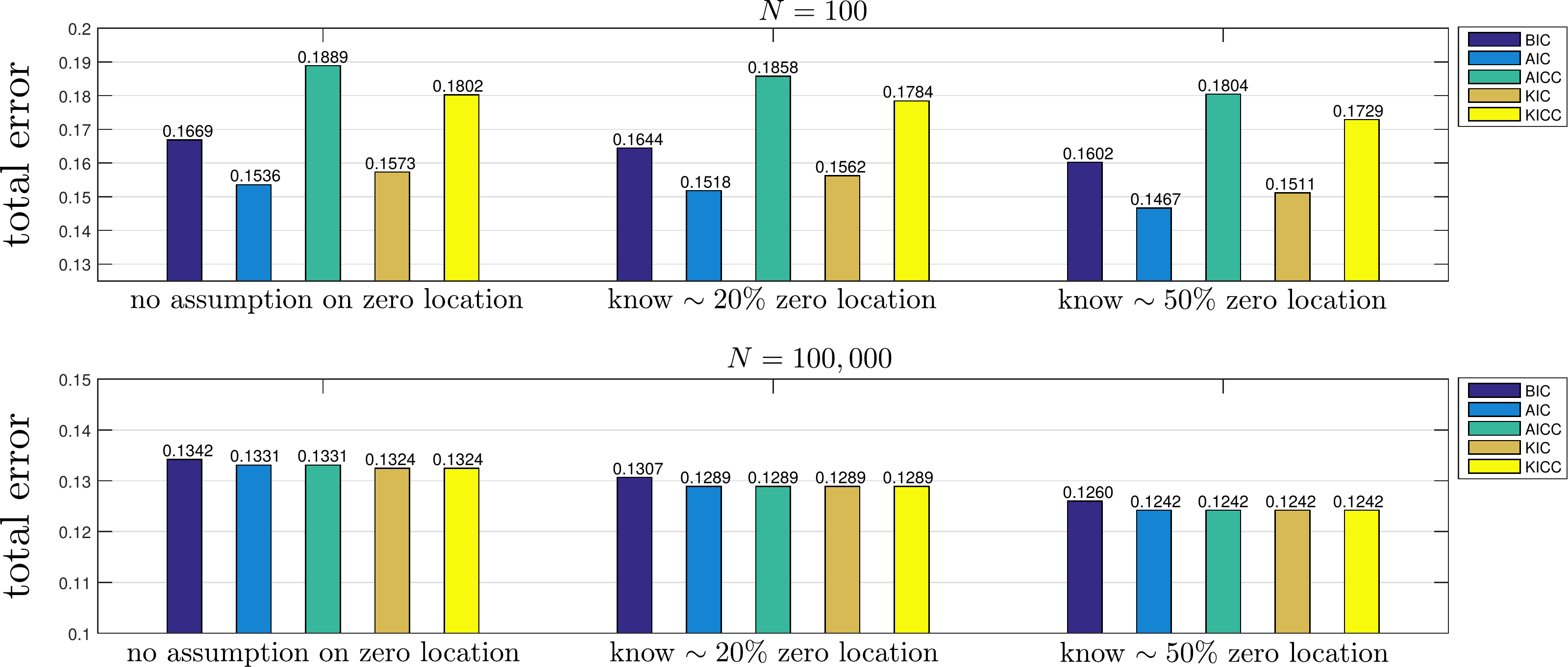}
            \caption{\textbf{Total error.} Main portion of the total error comes from FP so it tended to decrease with a greater assumption about the true zero location in $A_{\mathrm{true}}$.}
            \label{fig:total_error_Atrue_dense}
    \end{subfigure}

    	\caption{Averaged (a) FP, (b) FN and (c) total error from $50$ runs of a sample covariance matrix $S$, when $A_{\mathrm{true}}$ \textbf{was dense.} The AIC provided the minimum error when $N$ was small.}
	\label{fig:FP_FN_error_Adense}
\end{figure}

\begin{figure}
\centering
    \begin{subfigure}[b]{0.9\linewidth}         
            \includegraphics[width=\textwidth]{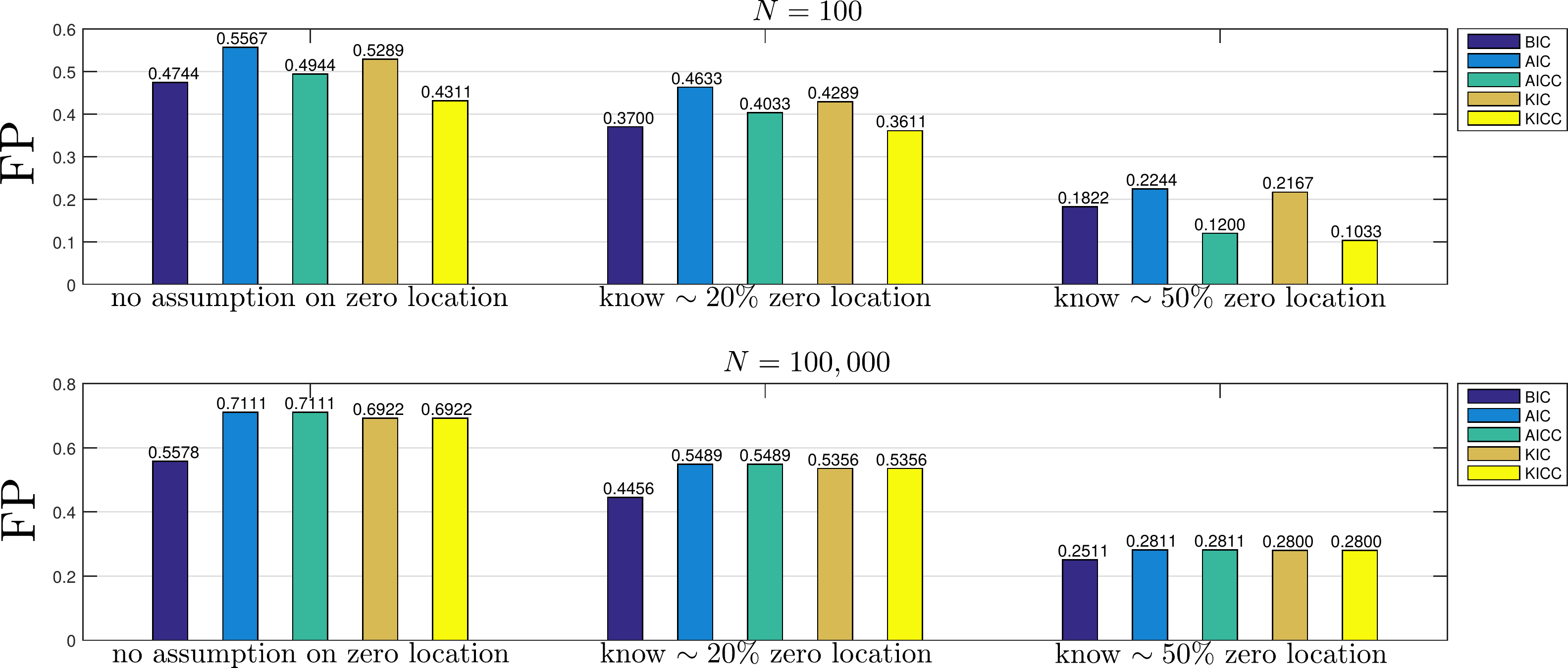}		
	\caption{\textbf{False Positive (FP) error.} The FP from all model selection criterions tended to decrease with more knowledge about the zero location in $A_{\mathrm{true}}$ in the estimation process, but increased when $N$ grew as all model selection criterions tended to select a denser $\hat{A}$. In the case of a small $N$, the BIC, AICc and KICc provided a better accuracy as the true model was sparse.}
		\label{fig:FP_Atrue_sparse}
	\end{subfigure}
	\\
    	\begin{subfigure}[b]{0.9\linewidth}
            \includegraphics[width=\textwidth]{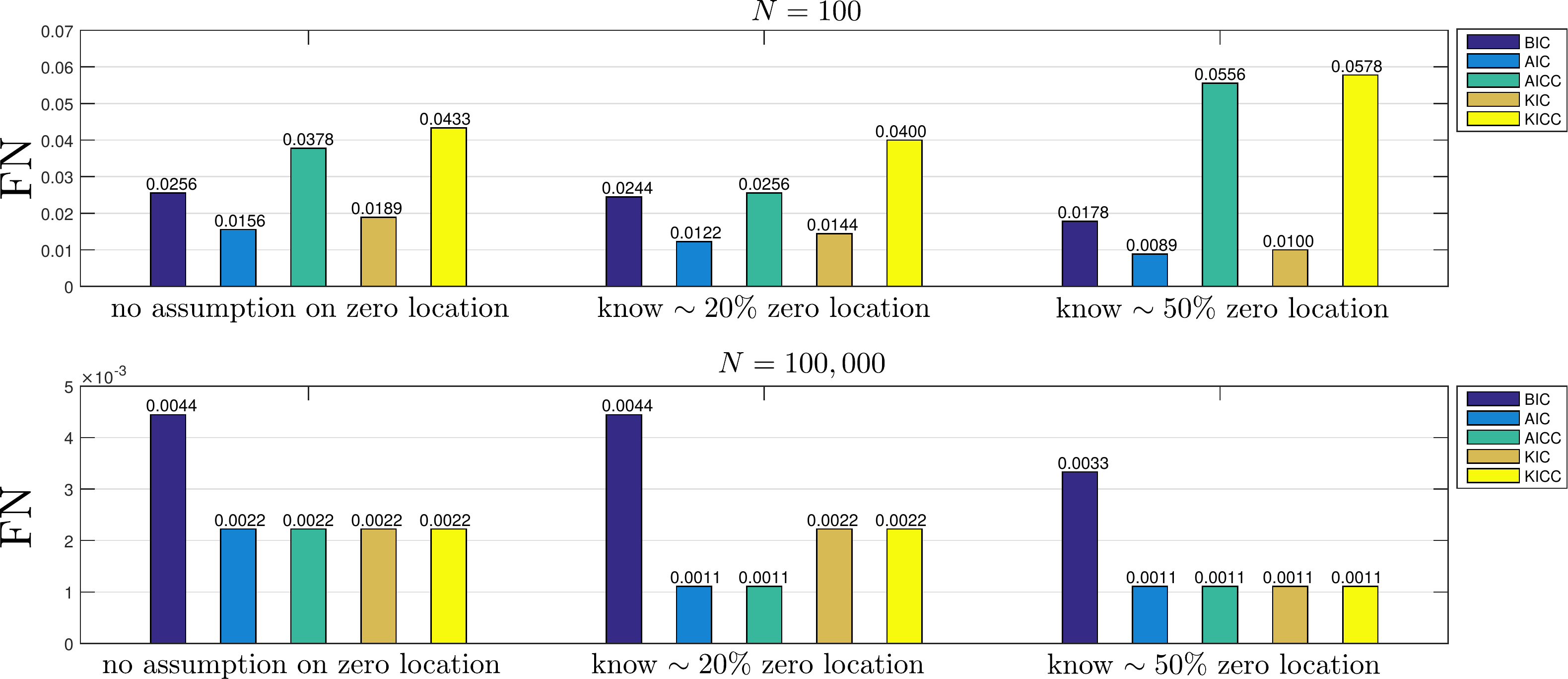}
		\caption{\textbf{False Negative (FN) error.} Using more knowledge about zero location in $A_{\mathrm{true}}$ in the estimation process barely affected the change of FN, but it improved when $N$ grew.}
		\label{fig:FN_Atrue_sparse}
	\end{subfigure}'
	\\
    	\begin{subfigure}[b]{0.9\linewidth}
            \includegraphics[width=\textwidth]{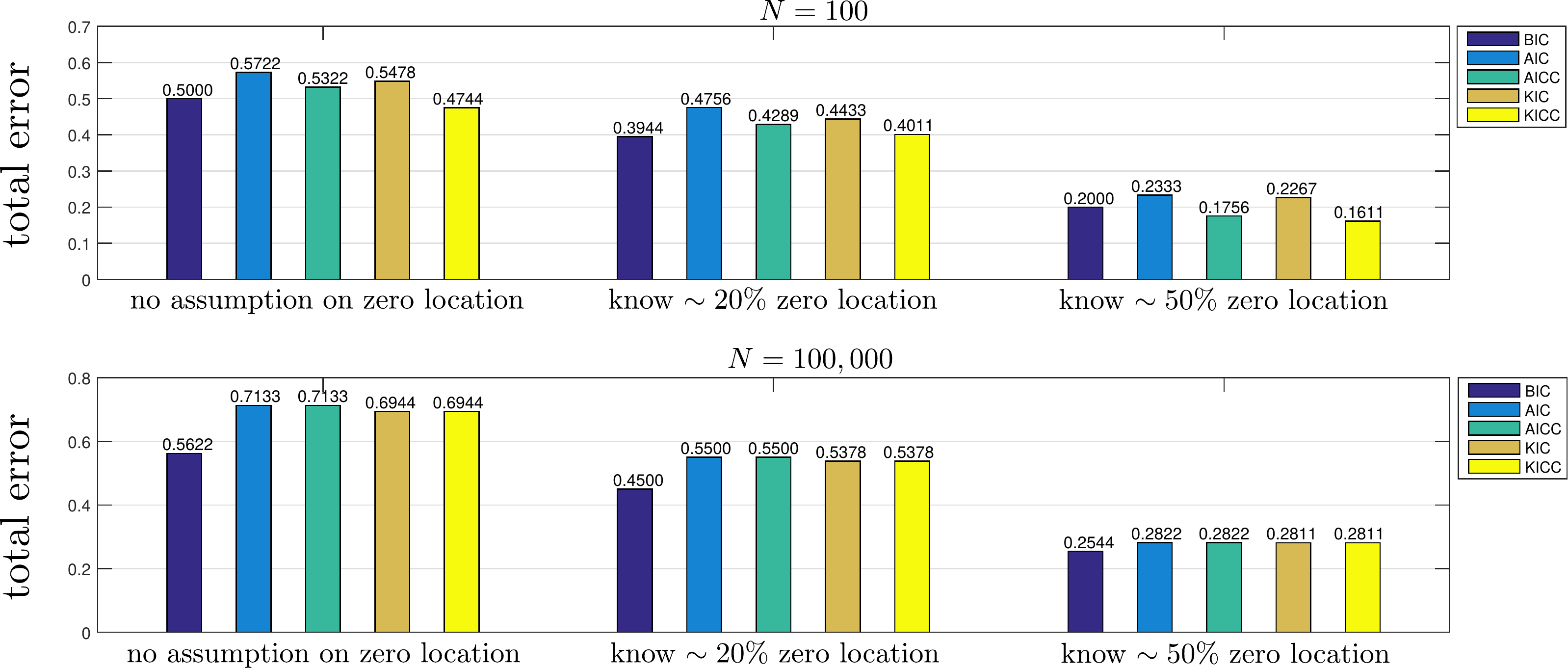}
		\caption{\textbf{Total error.} Main portion of the total error came from FP, which tended to decrease with more assumptions about the true zero location in $A_{\mathrm{true}}$.}
		\label{fig:total_error_Atrue_sparse}
	\end{subfigure}

	\caption{Averaged (a) FP, (b) FN and (c) total error from $50$ runs of sample covariance matrix $S$, when $A_{\mathrm{true}}$ \textbf{was sparse.} The BIC, AICc and KICc provided lower total errors when $N$ was small.}
	\label{fig:FP_FN_error_Asparse}
\end{figure}

\begin{table}[ht]
\centering
\caption{Chosen values of $\gamma$ by BIC scores, displayed as a scaling factor in $(0,1)$ of $\gamma_\mathrm{\max}$.}
\vspace{-2mm}
\begin{tabular}{|l| c|c|c|c|} \hline
\multirow{2}{*}{$\%$ assumed zeros in $A$}  & \multicolumn{2}{c|}{$N=100$} & \multicolumn{2}{|c|}{$N=100,000$} \\  \cline{2-5}
& dense $A_\mathrm{true}$ & sparse $A_\mathrm{true}$ & dense $A_\mathrm{true}$& sparse $A_\mathrm{true}$ \\ \hline
$0 \%$ of assumed zeros & $1.1 \times 10^{-3}$ & $2.62 \times 10^{-2}$ & $1.0 \times 10^{-4}$ & $ 6.9 \times 10^{-4}$ \\ 
$20 \%$ of assumed zeros & $1.4 \times 10^{-3}$ & $3.36 \times 10^{-2}$ & $1.0 \times 10^{-4}$ & $ 1.1 \times 10^{-3}$ \\ 
$50 \%$ of assumed zeros & $1.8 \times 10^{-3}$ & $3.36 \times 10^{-2}$ & $0.0 \times 10^{-4}$ & $ 1.1 \times 10^{-3}$ \\ \hline
\end{tabular}
\label{tab:BIC}
\end{table}

\begin{enumerate}
\item The sparsity density of $A_{\mathrm{true}}$. To this end, $A_{\mathrm{true}}$ was generated at the two sparsity levels of $50\%$ and $80\%$ and the relationship between the sparsity pattern of $\hat{A}$ that minimizes the BIC score and the error rate was observed. Comparing the FN from Figures~\ref{fig:FP_FN_error_Adense} and~\ref{fig:FP_FN_error_Asparse} for a moderate sample size ($N=100$) revealed that our method gave a lower FN when $A_{\mathrm{true}}$ was sparse and a lower FP when $A_{\mathrm{true}}$ was dense. Unavoidable errors as FP (when $A_{\mathrm{true}}$ is sparse) and FN (when $A_{\mathrm{true}}$ is dense) were commonly seen since these type of errors occur against the hypothesis of the true model. Moreover, when $A_{\mathrm{true}}$ was dense, using AIC leads to the minimum total error (Figure~\ref{fig:FP_FN_error_Adense}) since this score is prone to use a dense model (which agrees with the assumption of the true model). Similarly, when $A_{\mathrm{true}}$ was sparse, using the scores penalising more on the model complexity, and so BIC, AICc and KICc yielded a lower total error (Figure~\ref{fig:FP_FN_error_Asparse}). 

\item The number of samples. In the experiments, we use $N=100$ (moderate size) and $N=100,000$ (large sample size) to examine the asymptotic properties of the estimates. When $N$ was large, the selected $\gamma$ was closer to zero (Table~\ref{tab:BIC}) since the sparse SEM (as a regularized problem) should yield the solution closer to that of non-regularized problem. The selected $\gamma$ was also larger when the true model was sparse. Moreover, the FP was not improved (since the solutions are denser), but the FN obviously decreased when $N$ increased (Figures~\ref{fig:FP_FN_error_Adense} and~\ref{fig:FP_FN_error_Asparse}), showing that our regularized formulation was robust to FN errors.
\item The percentage of known zero locations in the estimation. To examine this factor, the experiments were performed with known zeros of $0\%, 20\%$ and $50\%$. The first two values corresponded to the problem with a negative df where the sparse SEM solution could be not unique, implying that the estimated zero pattern may not be as accurate as when knowing more zero locations. This, in principle, should affect the FP, as supported by Figures~\ref{fig:FP_FN_error_Adense} and~\ref{fig:FP_FN_error_Asparse}, where a greater knowledge about the true zero locations in $A_{\mathrm{true}}$ decreased the FP and overall total error, but FN seemed to be indifferent.

\item The choice of model selection scores. We considered the AIC (tended to choose dense models), AICc (adjusted for finite sample size), BIC, KIC and KICc scores (tended to choose simpler models). Figures~\ref{fig:FP_FN_error_Adense} and~\ref{fig:FP_FN_error_Asparse} supported that the AIC tended to yield the minimum total error when $A_{\mathrm{true}}$ is dense, and conversely, the BIC, AICc and KICc tended to provide lower total errors than those of the other criterions when $A_{\mathrm{true}}$ was sparse.
\end{enumerate}

\subsection{Comparison with an existing method}
\label{sec:regsem}
\begin{figure}[!ht]
\centerline{\includegraphics[width=0.5\columnwidth]{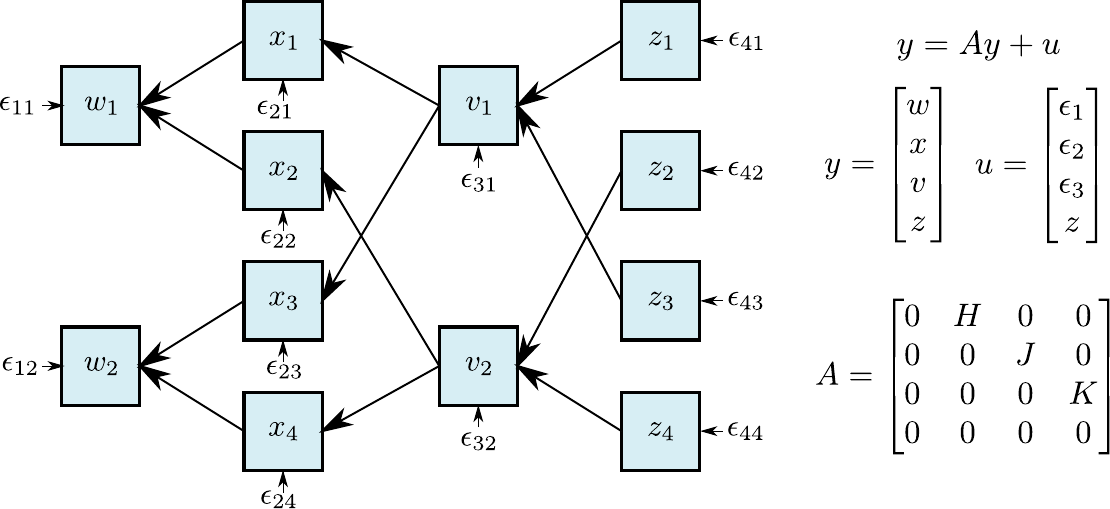}}
\caption{The recursive path model used in the comparison of our model with Regsem.} 
\label{fig:recursive_path}
\end{figure}

\texttt{Regsem} package in R was developed to solve sSEM with a regularization term, including both the ridge and the least absolute shrinkage and selection operator (lasso)~\cite{JGM16}. Regsem uses the RAM (Recticular Action Model) notation to derive an implied covariance matrix. The parameters of the general SEM were translated into three matrices: the filter matrix ($F$), the direct path matrix ($A$) and the covariance matrix of variables ($\Psi$). In detail, \texttt{regsem} can be used to solved an optimization problem shown in~\eqref{eq:regsem_lasso};
\begin{equation}
\begin{array}{ll}
\minimize & - \log \det \Sigma^{-1} + \Tr(S \Sigma^{-1}) - \log \det S - n + \gamma \sum\limits_{(i,j) \notin I_A} |A_{ij}| \\
\mbox{subject to} & \Sigma^{-1} = F^{-T} (I-A)^T \Psi^{-1} (I-A)F^{-1} , \quad P(A) = 0,
\end{array}
\label{eq:regsem_lasso}
\end{equation}
with variables $\Sigma$, $A$ and $\Psi$. It generalizes~\eqref{eq:sem_original} to include the (i) $\ell_1$-regularization penalty on $A$ and (ii) the filter matrix, $F$, since the RAM model also considers latent variables. However, our problem aimed to find the relationship among the observed variables only. Therefore, to be able to make a comparison between Regsem and our method, we explored the structure of a \emph{recursive path model}, which is the simplest model in RAM and is described by the set of linear equations: $w = Hx + \epsilon_1,\; x = Jv + \epsilon_2,\; v = Kz + \epsilon_3$,
where $w,x,v$ and $z$ are the observed vectors and $H,J$ and $K$ are the coefficient matrices. These equations can be written in $y=Ay +u$, as displayed in Figure~\ref{fig:recursive_path}, where $A$ has a certain sparse structure. In the experiment, we then generated $H_\mathrm{true},J_\mathrm{true}$ and $K_\mathrm{true}$ as a true description of the RAM model, which also implied $A_\mathrm{true}$. The measurements were generated from $y = (I-A_{\mathrm{true}})^{-1}u$, where $y \in \reals^{12}$ and $u \sim \mathcal{N}(0,1)$. 

\begin{figure}[ht]
\centering
    \begin{subfigure}[h]{0.45\linewidth}            
            \includegraphics[width=\linewidth]{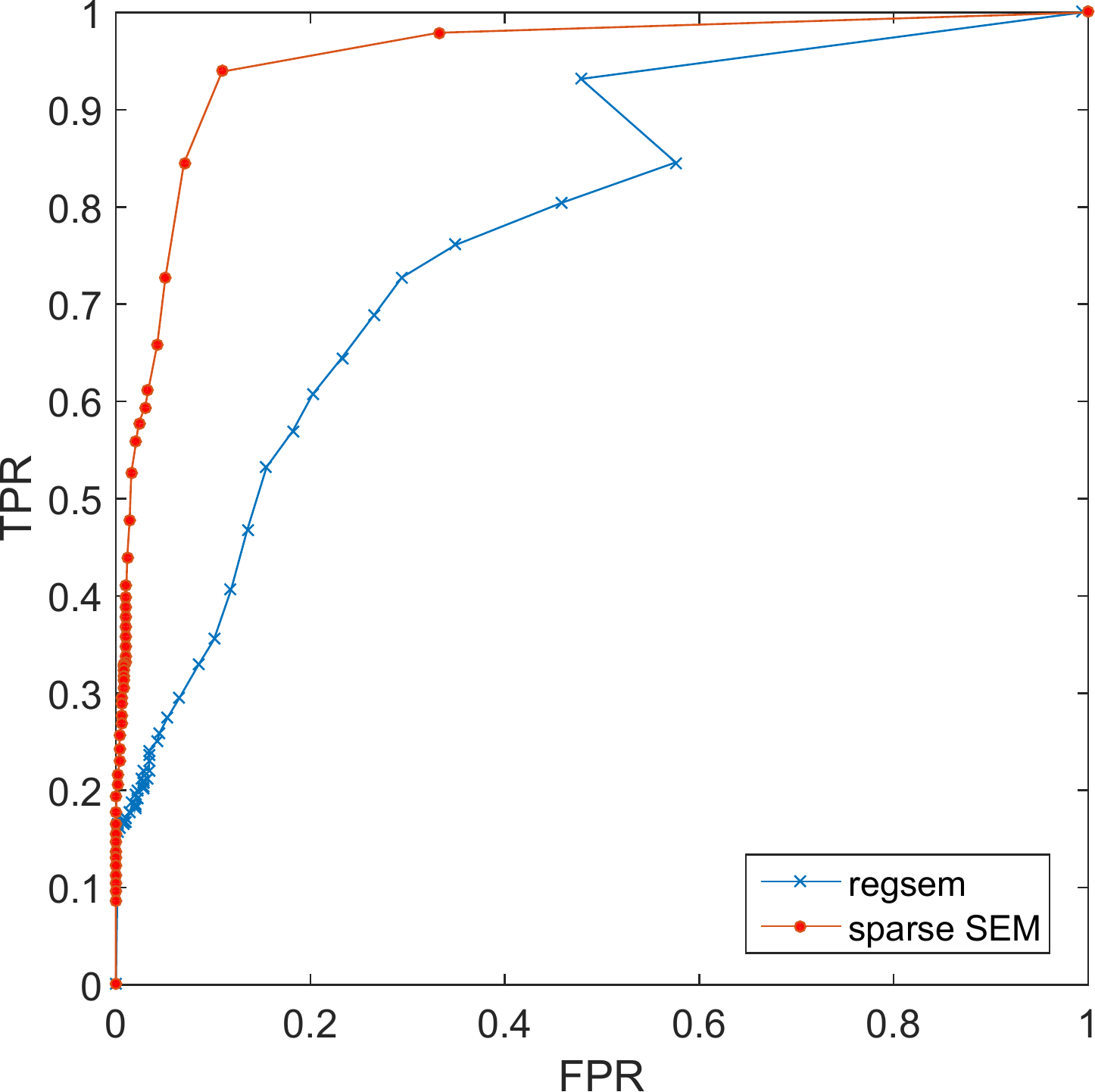}
            \caption{$N=100$.}
            \label{fig:compare_result_regsem_sparsesem_N100}
    \end{subfigure} \quad
    \begin{subfigure}[h]{0.45\linewidth}
            \includegraphics[width=\linewidth]{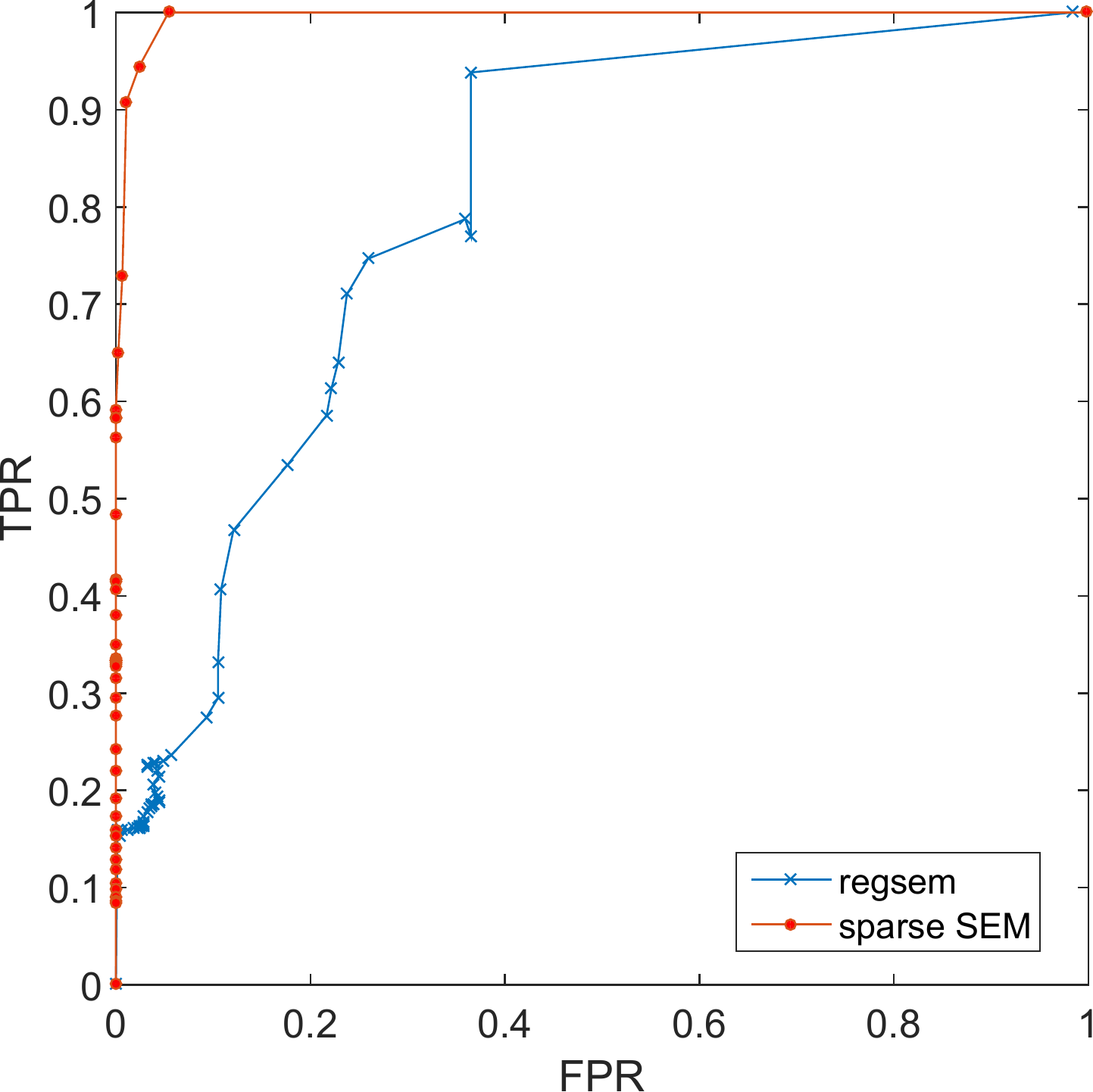}
            \caption{$N=1000$.}
            \label{fig:compare_result_regsem_sparsesem_N1000}
    \end{subfigure}

	\caption{The ROC curve from Regsem and sparse SEM averaged from 100 trials of a (a) moderate (N = 100) or (b) high (N = 1,000) sample size setting.}
	\label{fig:comparing_result_regsem_sparsesem}
\end{figure}

In our sparse SEM estimation process, the constraint $P(A) = 0$ was encoded according to the structure of $A_\mathrm{true}$ derived in Figure~\ref{fig:recursive_path} that corresponded to a df of $42$. For the Regsem, it was based on the \texttt{lavaan} package, which is a general SEM software program. This command required at least four arguments of our model generated by the \texttt{lavaan} library, type of estimation formulation (lasso or ridge regression), entries to be penalized and a penalty parameter. The argument $\text{pars}\_\text{pen}=c(1:24)$ in the \texttt{regsem} function was used to set the number of penalized entries in the path matrix $A$, which here was 24 entries, according to its structure in Figure~\ref{fig:recursive_path}. We also use this function to find a lower bound of the regularization parameter by selecting the minimum value of $\gamma$ that forced all entries in $A$ to be zero. For each data trial, we solved both the sparse SEM and regularized SEM in~\eqref{eq:regsem_lasso} by varying $50$ values of $\gamma$ in the range of $[0,\gamma_{\max}]$. 

The ROC curves averaged from $100$ trials (Figure~\ref{fig:comparing_result_regsem_sparsesem}) illustrate that our sparse SEM achieved a greater accuracy than Regsem in both the moderate and high sample size settings. Although the Regsem problem of~\eqref{eq:regsem_lasso} estimated the matrix variables in a more general model than our methods, when it comes to the special case that reduced to a path analysis problem and contained the observed variables only, our method, customized to a path analysis problem, should perform better. However, the Regsem model can be more advantageous when solving a more general SEM problem. 

\section{Real-world application}
\label{sec:realdata}

\subsection{Air pollution and weather data}
\label{sec:air}
We explored the relation structure among 11 climate variables (reenhouse gas (SO$_2$, NO$_2$, O$_3$ and CO), solar radiation, relative humidity, temperature, particulate matter (PM$_{10}$), pressure and wind speed) in residential sites of Bangkok, the capital of Thailand, during February 15, 2007 to May 15, 2009. The hourly data were measured from eight stations and were standardized to have a zero mean and unit variance. We split the data into training and validation sets at a $2:1$ ratio.

In the estimation process, if identifiability of the model was encouraged, zero df, as in~\eqref{eq:df}, should be obtained. Hence, it required setting $55$ entries in the estimated $A$ to be zero encoded in the constraint $P(X_2) = I$. In order to obtain a reasonable constraint, we performed a partial correlation analysis on $Y$ using \texttt{partialcorr} in \texttt{MATLAB} to find insignificant relations among the variables via the zero pattern in the sample partial correlation matrix (using a statistical test at a significance level of $0.01$) and imposed the corresponding entries in $A$ to zero. 

\begin{figure}[ht]
\centering
\includegraphics[width=0.8\linewidth]{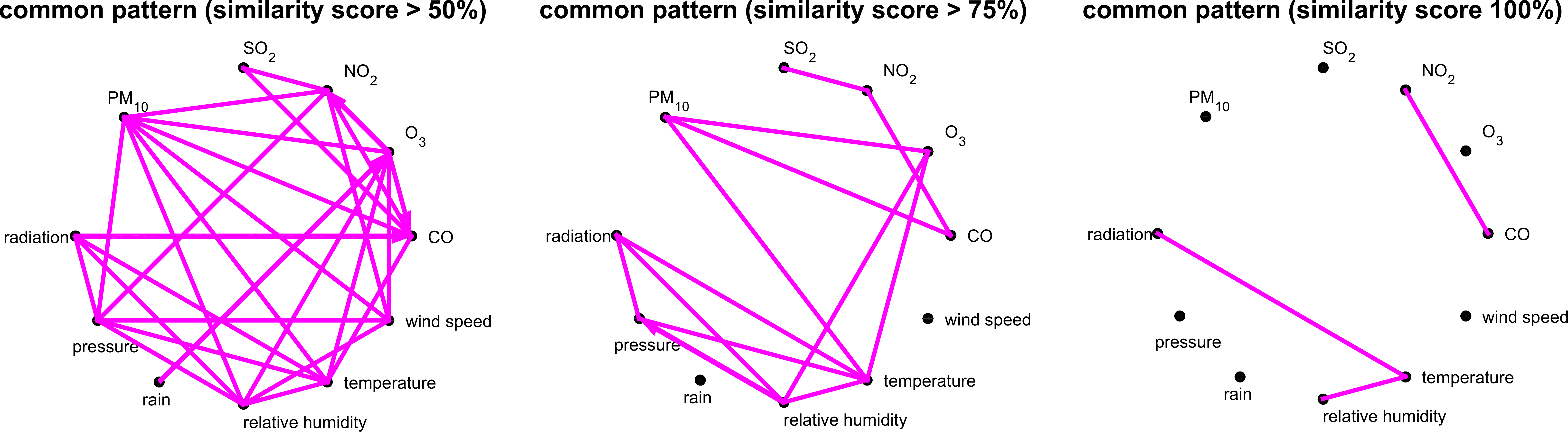}
	\caption{Common relationship pattern of variables from all stations with various similarity scores.}
	\label{fig:grap_commonnetwork}
\end{figure}

\begin{figure}
\centering
            \includegraphics[width=0.9\linewidth]{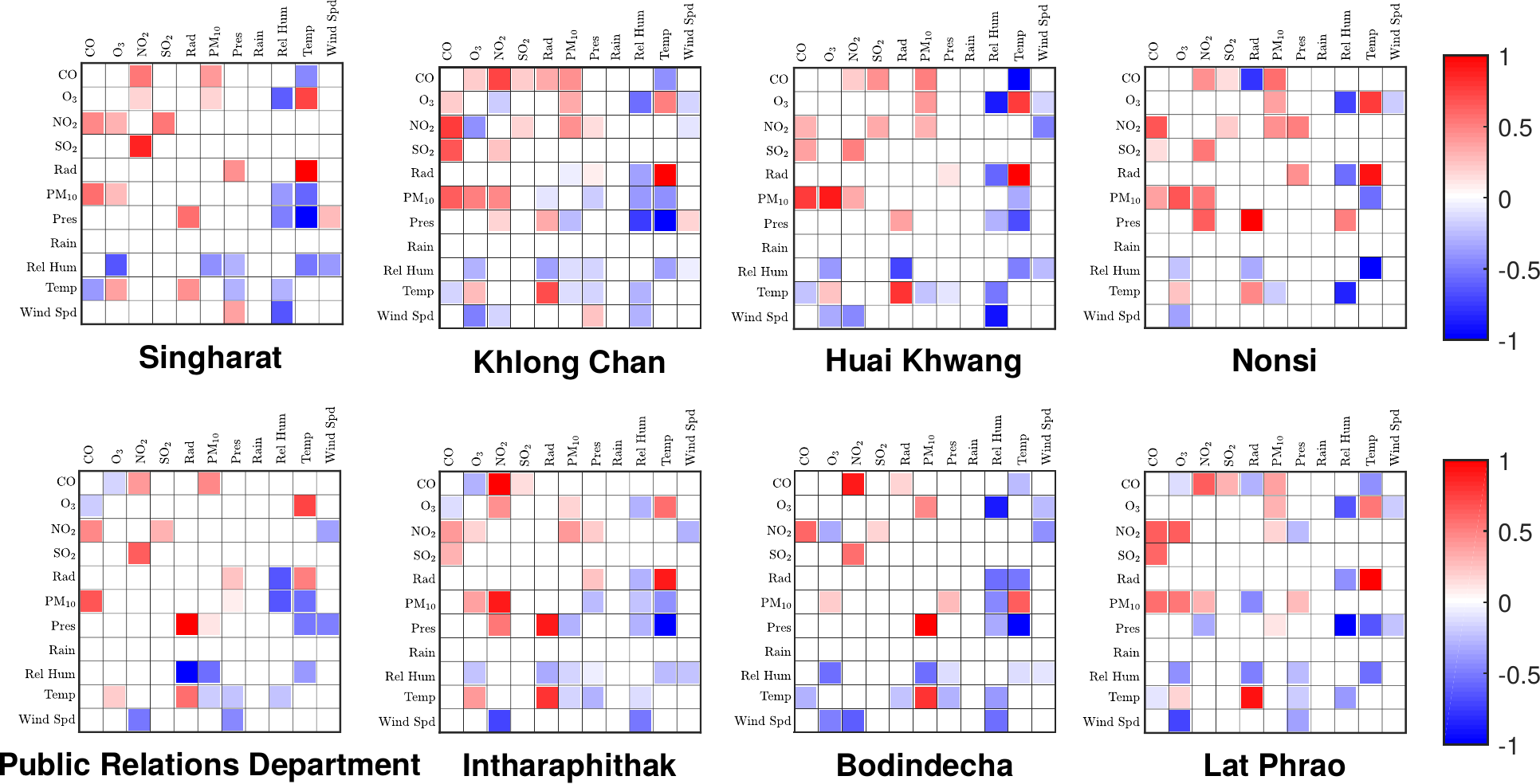}		
	\caption{The structure of the optimal path matrix from each station when zero constraints of $A$ from a partial correlation analysis is applied.}
	\label{fig:structure_assump_par_cor}
\end{figure}

The estimation process was performed as described in Section~\ref{sec:learning_path} and the path matrices from the model selection scheme of each station were selected. The nonzero entries in $A$ and their magnitude defined the graphical model and explained the relationships of variables and their strength of connections. These relationship patterns differed from station to station (Figure~\ref{fig:structure_assump_par_cor}) so we compared a common network from all stations using similarity scores ranging from $50\%$ to $100\%$ (Figure~\ref{fig:grap_commonnetwork}). The similarity score is the percentage of the number of common links from all stations relative to the number of all links. Dominant connections with a similarity score of $100\%$ were the:
\begin{center}
CO- NO$_2$, Radiation-Temperature, relative humidity (RH)-Temperature
\end{center}

Firstly, the strong relation between CO and NO$_2$ was supported by the known combustion reaction from car congestion in city areas. This appeared in our finding as the red color between the two variables from all stations in Figure~\ref{fig:structure_assump_par_cor}. Secondly, the relationship between the RH and temperature is known to be inverse to each other, given that the moisture content in the air is constant. As the temperature rises, the capacity of air to hold water increases, so if the actual amount of water in the air does not change, then the RH decreases. Mathematical expressions explaining the relationship between these two variables can be found in~\cite{Law:05}. Our result revealed an inverse relationship between the RH and temperature, as noted in the blue value from all stations in Figure~\ref{fig:structure_assump_par_cor}. Thirdly, the positive dependency between temperature and radiation, as indicated by the red color in Figure~\ref{fig:structure_assump_par_cor}, agrees with the natural characteristics of solar radiation, where as the sun rises up, the air temperature is increased. Nevertheless, in the area of meteorology, the daily solar radiation is typically explained by an increasing  diurnal air temperature range, $\Delta T$ (the difference between the maximum and minimum temperature)~\cite{LMLW+:09}.

In addition to the above relationships, the estimated structures with a similarity score of $75\%$ in Figure~\ref{fig:grap_commonnetwork} included the relations between
\begin{center}
CO-PM$_{10}$, PM$_{10}$- O$_3$, Temperature-O$_3$, RH-O$_3$, PM$_{10}$-Temperature, NO$_2$- SO$_2$, \\ RH-Radiation, RH $\rightarrow$ Pressure,
Radiation-Pressure and Temperature-Pressure.
\end{center}

The findings from~\cite{JaL:99} using a regression analysis explained that the wind speed and temperature had an inverse dependence on radiation attenuation, RH and pollution had a direct influence on radiation attenuation and that the reduction of solar radiation in the rainy season due to the pollution wash-out effect was not significantly different from the dry season. The relations among temperature, RH and radiation are consistent with several other researches \cite{JaL:99}, where the RH was found to have a direct influence on the radiation. 

The connection between PM$_{10}$ and O$_3$ corresponds to a previous finding~\cite{JCSH+:15}, although they considered PM$_{2.5}$ instead of PM$_{10}$. It is known that vehicle emissions are the main source of CO and NO$_2$, and diesel vehicles especially can emit particular matter. Moreover, areas containing burning process, such as fossil fuel combustion, can be major sources of NO$_2$ and SO$_2$. The airborne PM in Bangkok can be produced in areas of high temperature, such as from automobiles and biomass burning in residential areas~\cite{CNLK:08}. These facts relate to our findings of connections between CO-PM$_{10}$, PM$_{10}$-Temperature and NO$_2$- SO$_2$, where the last two relationships are present altogether in many stations (Figure~\ref{fig:structure_assump_par_cor}). As the temperature increases, the air expands and its density in that area is reduced, resulting in a decreasing variation in the air pressure. This agrees with the minus sign of the coefficient in the path matrix from temperature to pressure (Figure~\ref{fig:structure_assump_par_cor}). We conclude that the dependence structure learned from our model mostly agrees with known characteristics of environmental variables and with similar findings from previous studies.

\subsection{fMRI data}
In this experiment, we aimed to explore a common brain network between the control and autism groups under a resting-state condition learned from the Autism Brain Imaging Data Exchange data set~\cite{MYD+:14, abide}. This data set contains $1112$ fMRI images, collected across more than $24$ international brain imaging laboratories. We selected data trials that contain enough samples for estimation; $46$ images from the autism group and $40$ images from the control group, giving a total of $86$ images, provided by University of Michigan. 

The functional preprocessing was performed by the  Preprocessed Connectomes Project~\cite{pcp}, using a Configurable Pipeline for the Analysis of Connectomes. This performs a structural preprocessing of skull-stripping using \texttt{AFNI’s 3dSkullStrip}, three-issue type brain segmentation using \texttt{FAST} in \texttt{FSL} and skull-stripped brain normalisation to MNI152 with linear and non-linear registrations using \texttt{ANTs}. It then performs a slice timing correction and motion realignment, respectively. The image intensity was normalised by four-dimensional global mean and a band-pass filtering in the range of 0.01-0.1 Hz was applied. All images from every subject were transformed from the original to the MNI152 template. To reduce the data dimension, we averaged the time series over the region of interest (ROI) using an automated anatomical labeling template~\cite{TLPC+:02}, which was fractionated to a functional resolution $(3 \times 3 \times 3 \text {mm}^3)$ via nearest-neighbour interpolation.

\begin{figure}[!ht]
\centering
    \begin{subfigure}[b]{0.32\linewidth}            
            	\includegraphics[width=\textwidth]{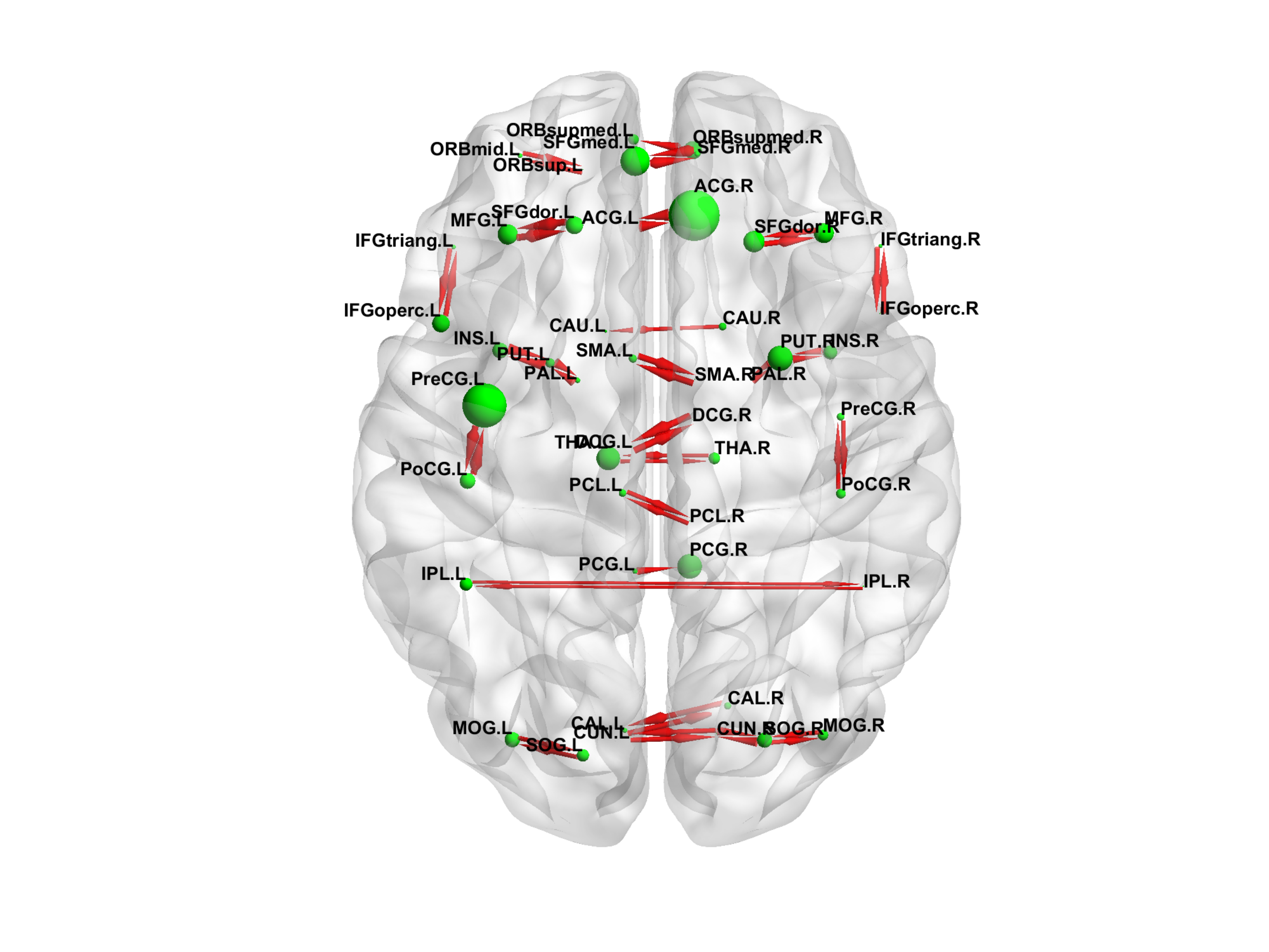}		
		\caption{Sparse network.}
		\label{fig:sparse_AU_3views}
	\end{subfigure}
    	\begin{subfigure}[b]{0.32\linewidth}
            	\includegraphics[width=\textwidth]{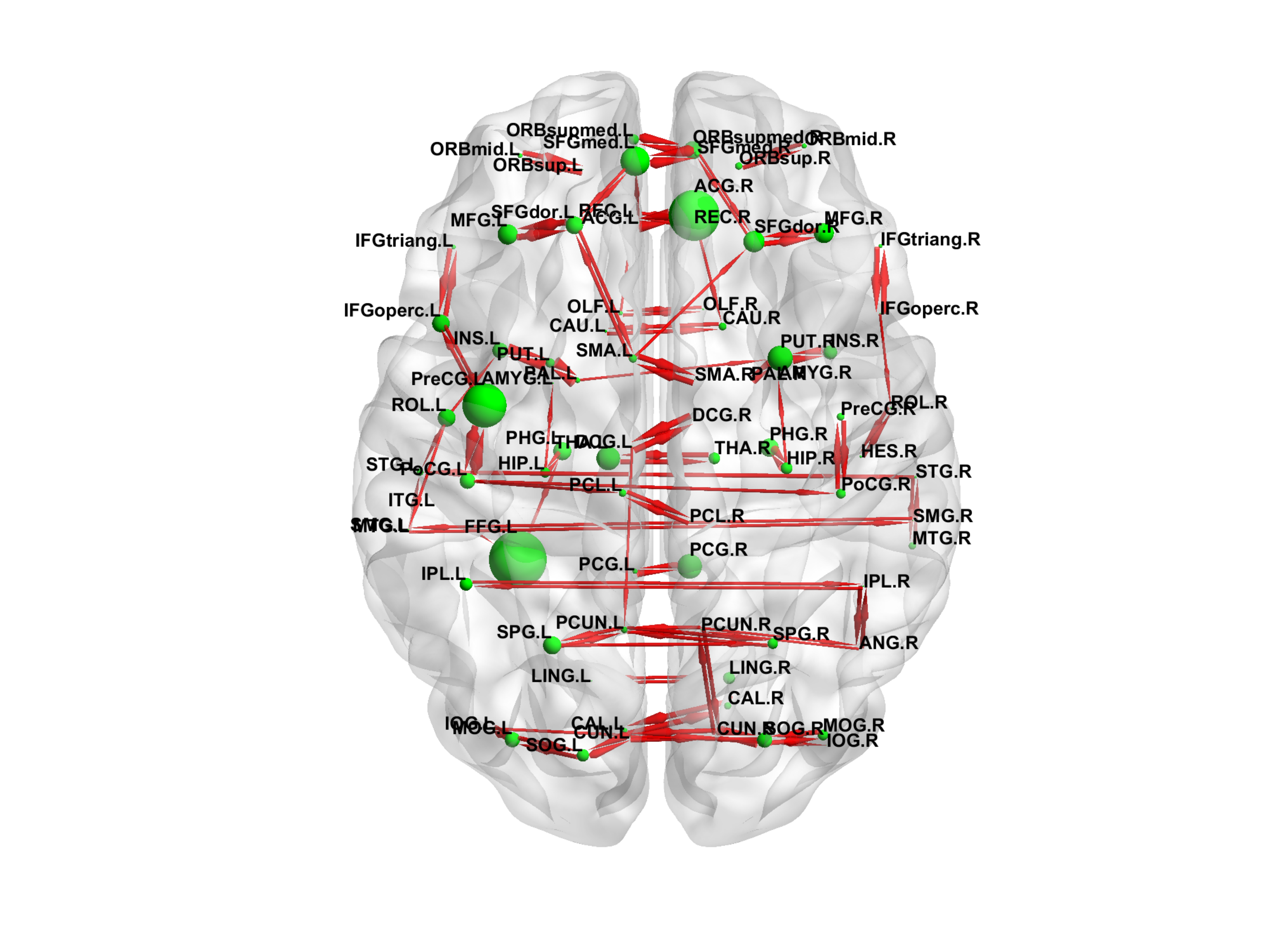}
		\caption{Moderately sparse network.}
		\label{fig:moderate_AU_3views}
	\end{subfigure}
    	\begin{subfigure}[b]{0.32\linewidth}
            	\includegraphics[width=\textwidth]{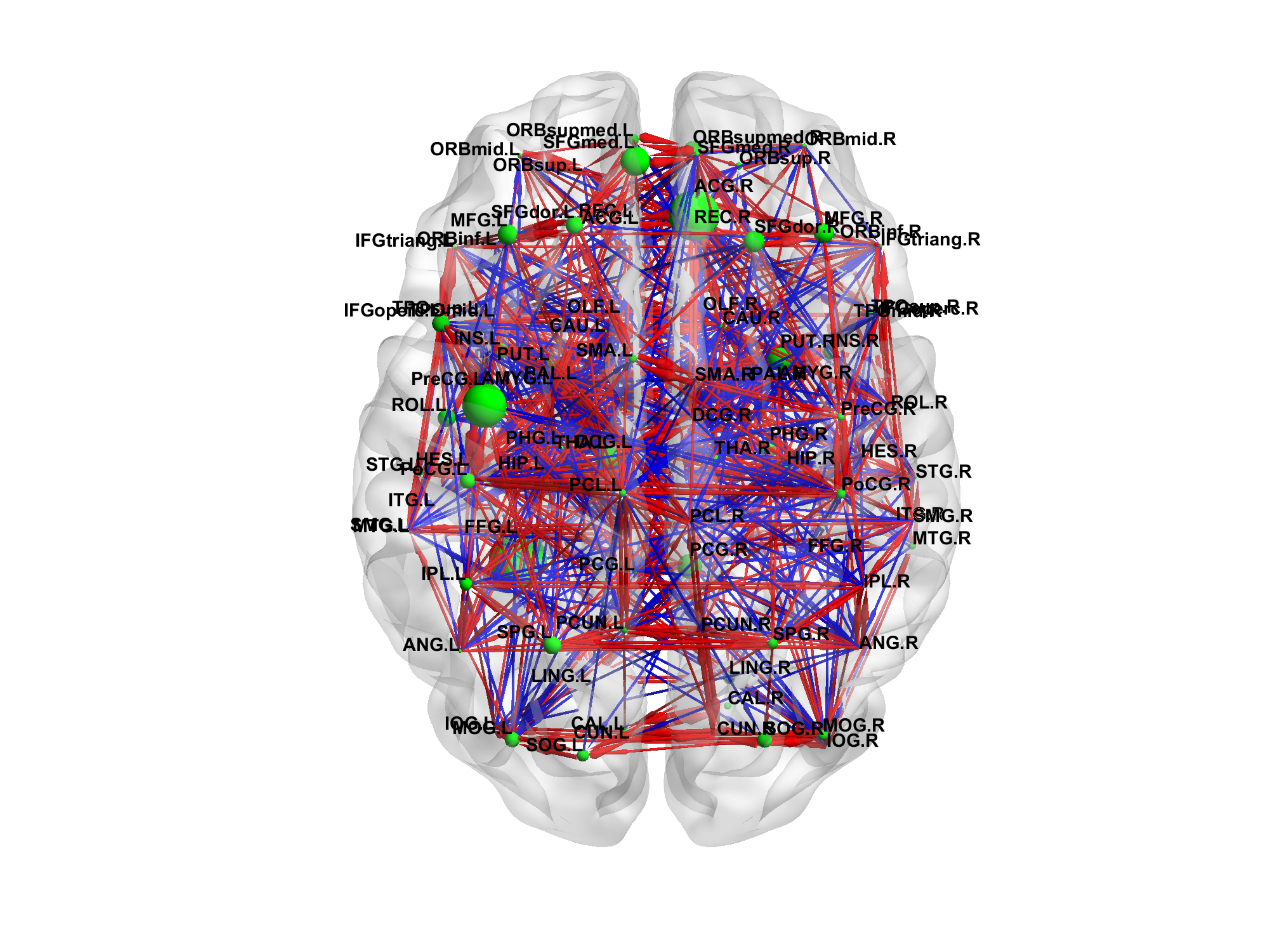}
		\caption{Dense network.}
		\label{fig:dense_AU_3viewsk}
	\end{subfigure}

	\caption{A common brain network from the \emph{autism group} at three graph density levels of (a) sparse, (b) moderately sparse and (c) dense. The link widths varied with the magnitudes of entries in the estimated path matrix. The red and blue link represent the positive and negative magnitude, respectively.}
	\label{fig:brain_network_AU}
\end{figure}

As a result, we obtained $Y \in \reals^{90 \times 249}$ from $90$ ROIs with $249$ time points. A prior assumption on the zero locations of $A$ was set as $P(A)=\diag(A) = 0$. Path matrices, $A \in \reals^{90 \times 90}$, as estimated from our approach, could be represented as a graphical model containing $90$ nodes whose labels were described in~\cite{FCZH+:15}. For each $\gamma$, we specified a common optimal path matrix by searching for positions of nonzero entry that appeared at a frequency of more than $90\%$ from all subjects in each group. The most frequently appeared links were denoted as significant connections between brain regions. We computed a common optimal path matrix by discarding insignificant entries and averaging only the significant entries over all subjects. We selected three values of $\gamma$, as $\gamma = 0.0025\gamma_{\max}, 0.0182\gamma_{\max}$ and $0.135\gamma_{\max}$, to produce three common path matrices with different density levels. Figure~\ref{fig:brain_network_AU} shows the brain networks from the autism groups, as constructed by \texttt{BrainNet Viewer}~\cite{BNV}. Our method generally provided nested graph structures as $\gamma$ varies, where noticeable connections in the sparse structure also existed in the moderately sparse and dense structures.

\begin{figure}[!hb]
\centering
    \begin{subfigure}[b]{0.47\linewidth}            
            	\includegraphics[width=\linewidth]{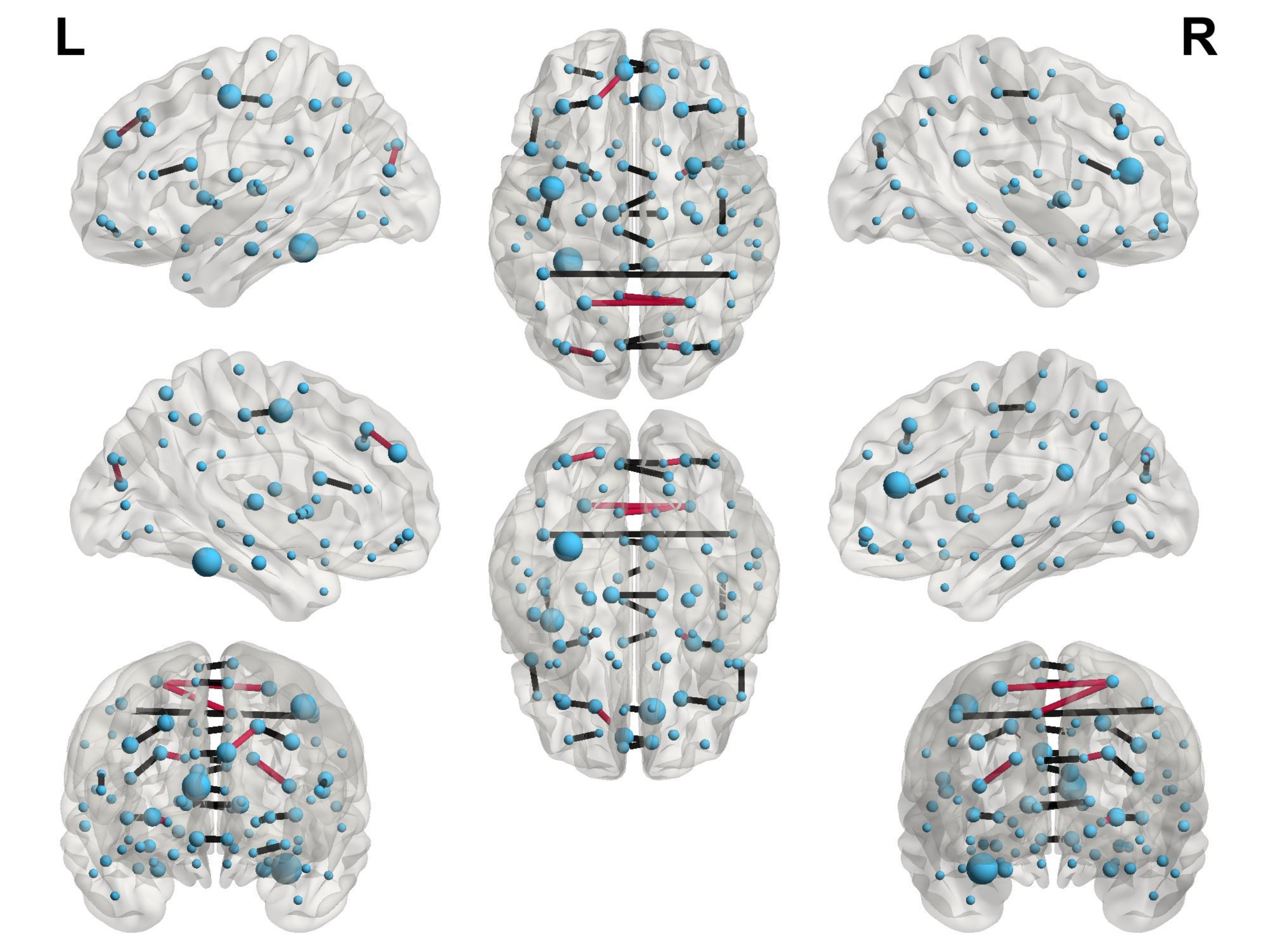}		
		\caption{A brain map of \textbf{autism} group}
		\label{fig:brainmap_AU}
	\end{subfigure}
    	\begin{subfigure}[b]{0.47\linewidth}       	
		\includegraphics[width=\linewidth]{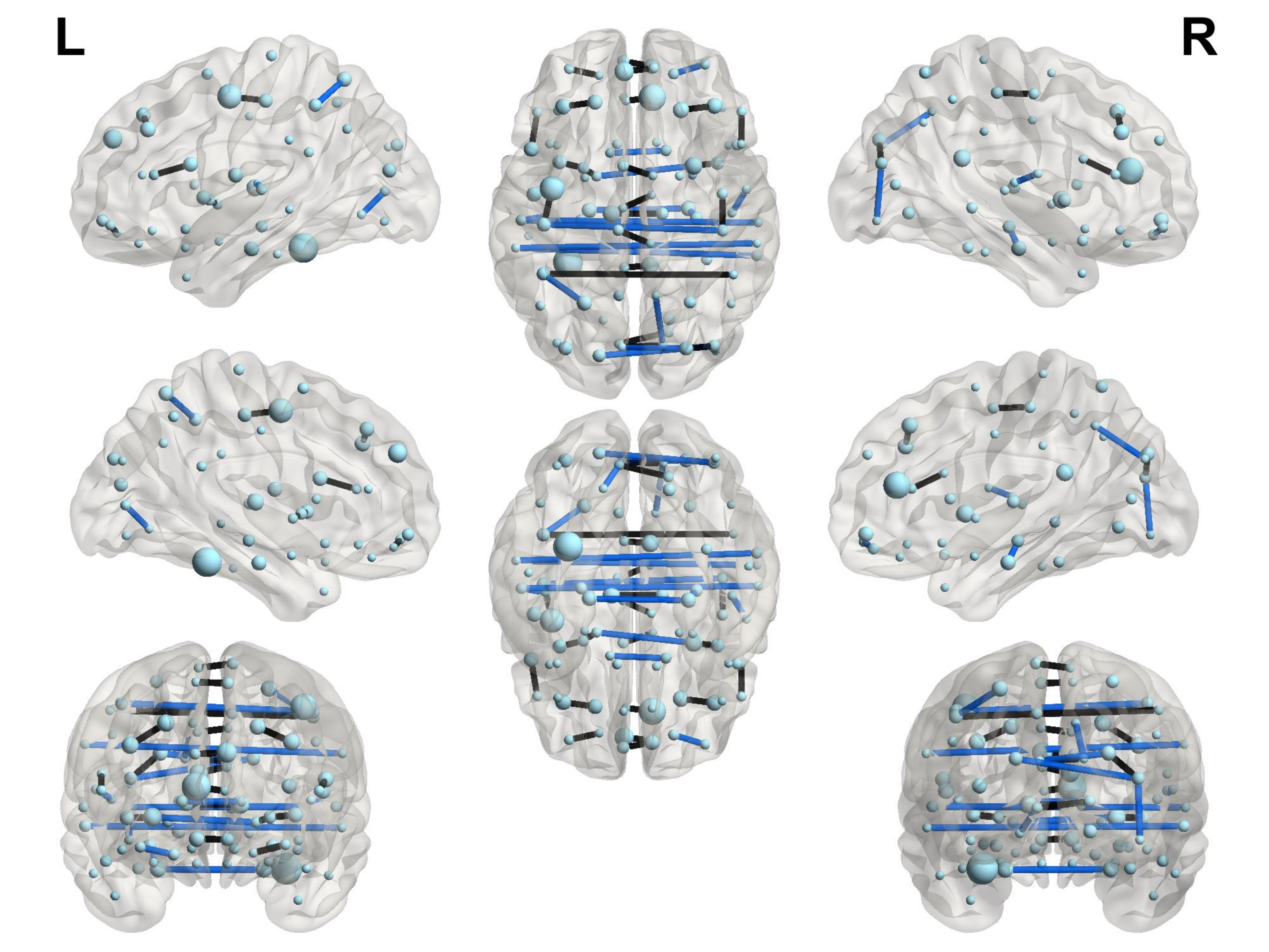}
		\caption{A brain map of \textbf{control} group}
		\label{fig:brainmap_CON}
	\end{subfigure}

	\caption{Comparing the brain maps between the autism and control groups. Black edges are common links that appear in both groups, red edges are ones found only in the autism group and blue edges are links detected only in the control group.}
	\label{fig:comparing_brainmap}
\end{figure}

To draw some conclusions on the network differences between the two groups, we selected the sparsest structures and compared the brain graphs in Figure~\ref{fig:comparing_brainmap}. An interesting observation was that the autism brain graph has a lower number of connections. The connections found only in the \emph{control} group were: \\

MFG(R) $\leftarrow$ MFG(L), ORBsupmed(R) $\leftarrow$ ORBsupmed(L), DCG(R) $\leftarrow$ ACG(R), PHG(L) $\leftarrow$ PHG(R), HIP(R) $\leftarrow$ PHG(R), PCUN(R) $\leftarrow$ CUN(R), PCUN(L) $\leftarrow$ SPG(L), IPL(R) $\leftarrow$ SPG(R), CUN(R) $\leftarrow$ PCUN(R), ROL(R) $\leftarrow$ HES(R), PUT(R) $\leftarrow$ PAL(L), STG(R) $\leftarrow$ HES(R), TPOmid(R) $\leftarrow$ TPOmid(L), ITG(R) $\leftarrow$ ITG(L), ORBmid(R) $\leftrightarrow$ ORBsup(R), LING(L) $\leftrightarrow$ CAL(L), SOG(R) $\leftrightarrow$ CUN(L), MOG(R) $\leftrightarrow$ SOG(L), IOG(R) $\leftrightarrow$ MOG(R), PoCG(L) $\leftrightarrow$ PoCG(R), IPL(L) $\leftrightarrow$ SPG(L), SMG(L) $\leftrightarrow$ SMG(R), CAU(L) $\leftrightarrow$ CAU(R), STG(L) $\leftrightarrow$ HES(L), STG(L) $\leftrightarrow$ STG(R) and MTG(L) $\leftrightarrow$ MTG(R).

From the list above (see abbreviations in~\cite{FCZH+:15}), it is worth comparing our results with previous studies on some of the missing connections that are no longer present in the autism group, which are among the temporal gyrus group (ITG(R) $\leftarrow$ ITG(L), STG(L) $\leftrightarrow$ STG(R) and MTG(L) $\leftrightarrow$ MTG(R)).  The inferior temporal gyrus (ITG) contributions are involved in processing of visual stimuli in object recognition~\cite{KW:12}, while the superior temporal gyrus (STG) is involved in the auditory process and developing language, but has recently also been implicated as a key factor in social cognition~\cite{BMNO+:07}, an important consideration of autism patients. The main function of the middle temporal gyrus (MTG) is not exactly known, but it helps in some processes about recognition of known faces and accessing word meaning while reading~\cite{AH:13}. These results are supported by  experimental results~\cite{CRZS+:17} in which a decreasing functional connectivity among these areas was observed in the autism group. These experiments~\cite{CRZS+:17} were designed to study the disruptive change in the state and strength of connectivities in the autism group, referred to as abnormal connectivities, which were also found in ITG and STG area as well. This is consistent with a previous report~\cite{CRZJ+:15}, which concluded that the MTG region is implicatexd in facial expressions, gesture representation impairments and theory of mind impairments in autism. We also found that relations from the precuneus (PCUN), the basal ganglia and the anterior cingulate cortex [PCUN(R) $\leftrightarrow$ CUN(R), CAU(L) $\leftrightarrow$ CAU(R) and DCG(R) $\leftarrow$  ACG(R)] were still missing from the brain network in the autism group. Caudate nucleus (CAU), one of the components in the basal ganglia, affects many nonmotor functions, such as procedural learning~\cite{NHM:09}, while the PCUN  is involved in self-consciousness, such as reflective self-awareness~\cite{KNL:02}. These results agree with~\cite{SWQ:17}, where brain networks were expressed as conditional independence graphs and revealed that, in the autism group, the edges linking to  the precuneus, the basal ganglia, the anterior cingulate cortex and the medial frontal cortex were mostly affected. Moreover, underconnectivity among the ROIs linking to the PCUN were reported in the autism group~\cite{CKKJ:06}.

In conclusion, our result represents that within the brain network some circuits relating to cognitive processes, such as social  interaction, face and image recognition, learning process or working memory, are missing from the autism group but exist in the control group. Our findings are similar to many previous studies.

\section{Conclusions}
\label{sec:conclusion}
This paper proposed two convex formulations for confirmatory and exploratory SEM, which have applications in learning causality among variables based on path analysis models. Our formulations relaxed the quadratic matrix equality into an inequality, where the solution to this formulation can be a solution to the original problem under homoskedastic assumptions of noise in the model. The proposed scheme of exploratory SEM exploits the feature of sparse estimation, introduced by adding an $\ell_1$ regularization to the estimation function. Causality structures encoded from sparsity patterns of the path matrix can be obtained by sweeping values of regularization parameters in a specific range, which was derived analytically. Applying efficient efficient proximal parallel algorithm allowed the problems to be solved at a large-scale, including hundreds of variables, while it is not common to see existing SEM softwares solve such cases. Numerical results showed that (i) the percentage of known zeros of $A$ was a key factor for decreasing the FP, regardless of the sparsity density in the true model, (ii) FN could be improved the most by increasing the sample size, (iii) the total error mainly came from FP and (iv) the choice of model selection criterions depended on the sparsity density in the true model. That is, BIC, AICc and KICc provided a better accuracy when the true path matrix was sparse, while AIC performed best when the true model was dense. Comparison between our approach with Regsem method was set in a fair setting, while considering a special case of the SEM model. Our method yielded higher accuracies where the desirable results could have benefited from either our formulation or the algorithms. Results from real data sets showed that the causality structures findings coincided with previous studies of applications. Most relations between climate variables learned from our model can be supported from environmental facts and previous research findings in the literature. While in the brain studies, although a ground-truth causal network of brain activities is not completely known, our findings were comparable with the literature with some extents. The brain networks from the control and autism groups are different in some brain regions, including the temporal gyrus group, PCUN and CAU areas.

\section{Acknowledgment}
We thank the Department of City Planning, Bangkok Metropolitan Administration and Pollution Control Department, Bangkok, for providing the relevant information on air pollution. This researh project was financially supported by a Chula Engineering research grant.
\bibliographystyle{alpha}
\bibliography{books,cvx-refs2,cvx-refs3,sem_ref,sparse,fmri_ref,jss}

\section{Mathematical proofs}
\subsection{Proof of Proposition~\ref{prop:alphac}}
\label{proof:alphac}
We applied a generalization of Farka's lemma to semidefinite programming~\cite{BoV:04}.
\begin{lemma}~\cite{BoV:04} The system: $U \succeq 0, \quad \Tr(GU) > 0, \quad \Tr(F_iU) = 0, i=1,2,\ldots,n$, is a strong alternative for the nonstrict linear matrix inequality (LMI): $\sum_{i=1}^n x_i F_i + G \preceq 0 $, if the matrices $F_i$ satisfy $\sum_{i=1}^n v_i F_i \succeq 0 $ implies that $\sum_{i=1}^n v_i F_i =0$.
\label{lem:farka}
\end{lemma}
Using the notation of $\tilde{A} = I-\tilde{X}_2$, the feasibility problem of~\eqref{eq:invS_feas} was expressed as an LMI as shown in~\eqref{eq:zerosol_block};
\begin{equation}
\begin{bmatrix} S^{-1} & (I-\tilde{A})^T & 0 \\ I-\tilde{A} & X_4 & 0 \\ 0 & 0 & \alpha I - X_4 \end{bmatrix} \succeq 0 \Longleftrightarrow \underbrace{\begin{bmatrix} -S^{-1} & -I & 0\\ -I & 0 & 0 \\ 0 & 0 & -\alpha I \end{bmatrix}}_{G} +  \underbrace{ \begin{bmatrix} 0 & \tilde{A}^T & 0 \\ \tilde{A} & 0 & 0 \\ 0 & 0 & 0 \end{bmatrix}}_{\sum_{ij} \tilde{A}_{ij} F_{ij}} +   \underbrace{\begin{bmatrix} 0 & 0 & 0 \\ 0 & -X_4 & 0 \\ 0 & 0 & X_4 \end{bmatrix}}_{\sum_{ij} (X_4)_{ij} H_{ij}} \preceq 0.
\label{eq:zerosol_block}
\end{equation}
The matrices $F_{ij}$ and $H_{ij}$ are a common choice of standard basis matrices that make up the above summation. In detail, let $E_{ij}$ be a standard basis matrix for the set of $n \times n$ matrices with zero diagonals and $S_{ij}$ be a standard basis matrix for $\symm^n$. In other words, the entries of $E_{ij}$ are all zero except that the $(i,j)$ entry is $1$. Similarly, the entries of $S_{ij}$ are all zero except that $(i,j)$ and $(j,i)$ entries are $1$. Note that $\tilde{A}$ only contains non-zero entries for $(i,j)\notin I_A$, so the expressions of $F_{ij}$ and $H_{ij}$ are 
\[
F_{ij} = \begin{bmatrix} 0 & E_{ij}^T & 0 \\ E_{ij} & 0 & 0 \\ 0 & 0 & 0 \end{bmatrix}, \;\text{for} \;\;(i,j) \notin I_A,\;\;
H_{ij} = \begin{bmatrix} 0 & 0 & 0 \\ 0 & -S_{ij} & 0 \\ 0 & 0 & S_{ij} \end{bmatrix},\;\text{for}\;\; i \geq j=1,2,\ldots,n.
\]
From Lemma~\ref{lem:farka}, the LMI of~\eqref{eq:zerosol_block} has no solution if and only if $\exists U \succeq 0, U \neq 0$ such that $\Tr(GU) > 0, \quad  \Tr(F_{ij} U ) = 0, \;\;\text{for $(i,j)\notin I_A$},\quad \Tr(H_{ij}U) = 0, \;\;\text{for $i \geq j$}$. In what follows, we show that there always exists such a matrix $U$ under the condition $\alpha \leq \alpha_c$. For scalars $\gamma$ and $\beta$ with $\beta \geq 0$ and $\gamma \neq 0$, we construct $U \succ 0$ of the form:
\[
U = \begin{bmatrix} (\gamma^2/\beta)I & \gamma I & 0 \\ \gamma I & \beta I & 0 \\ 0 & 0 & \beta I\end{bmatrix}.
\]
With this choice, we can easily check that $\Tr(F_{ij}U)= 0$ regardless of the choice of $I_A$ (as long as $I_A$ contains the indices of diagonal entries of $A$), and that $\Tr(H_{ij} U) = 0$. Moreover, the condition $\Tr(GU) > 0$ is expressed as in~\eqref{eq:poly_gamma},
\begin{equation}
\frac{\Tr(S^{-1})}{\beta} \left ( \gamma^2  + \frac{2n\beta}{\Tr(S^{-1})} \cdot \gamma + \frac{n}{\Tr(S^{-1})} \alpha \beta^2 \right ) < 0.
\label{eq:poly_gamma}
\end{equation}
The above quadratic polynomial in $\gamma$ can be expressed in terms of $\alpha$ and $\alpha_c$ as $\gamma^2 + 2 \alpha_c \beta \gamma + \alpha \alpha_c \beta^2 \leq 0$. Therefore, if $\alpha \leq \alpha_c$ then we can always choose any negative real value of $\gamma $ in the interval $\left ( -\alpha_c \beta (1+ \sqrt{1-\alpha/\alpha_c}) , -\alpha_c \beta (1- \sqrt{1-\alpha/\alpha_c}) \right )$, so that~\eqref{eq:poly_gamma} is satisfied with strict inequality. Lastly, we also see that $\sum_{ij} \tilde{A}_{ij} F_{ij} + \sum_{ij} (X_4)_{ij} H_{ij}  \succeq 0$ only implies that $X_4 = 0$ and $\tilde{A} = 0$ because its leading $(1,1)$ block is zero. From Farka's lemma, this concludes that if $\alpha \leq \alpha_c$ the feasibility problem of~\eqref{eq:zerosol_block} always has no solution. 
\subsection{Derivation of $\gamma_{\max}$}
\label{app:gammamax}
We will show that there exists a critical value of $\gamma$, denoted by $\gamma_\mathrm{\max}$, such that if $\gamma \geq \gamma_\mathrm{\max}$, then the optimal solution of $A$ obtained by $A=I-X_2$ in~\eqref{eq:seml1_primal} is zero. The derivation of $\gamma_\mathrm{\max}$ is, in fact, derived from the KKT conditions of~\eqref{eq:seml1_primal} and under an \emph{assumption} that the optimal primal solution is \emph{low rank}.
\paragraph{KKT conditions of~\eqref{eq:seml1_primal}.} If strong duality holds, $X$ and $Z$ are optimal if and only if the following conditions hold. 
\begin{itemize}
\item \textbf{Zero gradient of the Lagrangian:} $X_1 = (S-Z_1)^{-1}$.
\item \textbf{Primal and Dual feasibility:} $ X \succeq 0$, \;\; $ 0 \preceq X_4 \preceq \alpha I$, \;\; $Z  \succeq  0, \;\;\Vert P^c(Z_2) \Vert_\infty \leq \gamma$.
\item \textbf{Complementary slackness condition:} $ZX=0$ and $Z_4(X_4-\alpha I) = 0$.
\end{itemize}
If the optimal $X$ has rank $n$, then it follows from the complementary slackness condition that $\Rank Z = n$ and $\Rank Z_4 = n$, so $Z_4$ is invertible. This further implies from the slackness condition: $Z_4(X_4-\alpha I)=0$ that $X_4= \alpha I$. Since we aim to characterize the dual feasibility condition when we obtain the sparsest solution of $A$, we set $A=0$ (or equivalently $X_2 = I$) in the optimal solution, then $X = \begin{bmatrix} X_1 & I \\ I & \alpha I \end{bmatrix}$. We see that $X$ has rank $n$ if and only if $X_1 = (1/\alpha)I$. The zero gradient of the Lagrangian condition, $Z_1 = S-X_1^{-1} = S-\alpha I$, is substituted in the slackness condition: $ZX =0$,
\[
 \begin{bmatrix} S-\alpha I & Z_2^T \\ Z_2 & Z_4 \end{bmatrix} \begin{bmatrix} X_1 & I \\ I & \alpha I \end{bmatrix} = 0
\]
from which we can solve for $Z_2$ as $Z_2 = (1/\alpha) ( \alpha I - S)$ and the dual feasibility condition becomes $\gamma \geq (1/\alpha) \Vert P^c( \alpha I - S) \Vert_\infty$. In conclusion, we showed that if $X_2=I$ is the optimal solution to~\eqref{eq:seml1_primal} and the optimal $X$ has rank $n$, then $\gamma$ must exceed a particular value. The KKT conditions are sufficient and necessary conditions for the optimality of a convex problem. As a result, we can set $\gamma_\mathrm{\max} = (1/\alpha) \Vert P^c( \alpha I - S) \Vert_\infty$ as the critical value of $\gamma$, and conclude that for any $\gamma \geq \gamma_\mathrm{\max}$, the optimal solution $A$ must be zero.


\subsection{Solution of a scaled sparse SEM}
\label{sec:scaled_sem}
We provide a proof of Proposition~\ref{prop:scaled_sem}. We show that if $(X,Z)$ satisfies the KKT conditions of the unscaled problem in~\eqref{eq:seml1_primal} using parameter $(S,\alpha)$ then $(\tilde{X}, \tilde{Z})$ provided in the statement also satisfies the KKT conditions of the scaled problem~\eqref{eq:seml1_primal} using parameter $(\tilde{S},\tilde{\alpha})$. 
\paragraph{Primal feasibility.} Given that $X \succeq 0$ which is equivalent to $ X_1 - X_2^T X_4^{-1} X_2 \succeq 0 ,\quad X_4 \succ 0$, by the Schur complement, then for any $\beta > 0$, it follows that the previous inequalities are preserved: $\beta X_4 \succ 0 ,\quad X_1/\beta - X_2^T (\beta X_4)^{-1} X_2) \succeq 0$ which are equivalent to $\tilde{X} \succ 0$ by the Schur complement when $\tilde{X}$ is given in~\eqref{eq:scaled_XZ}. Moreover, it is obvious that if $X_4 \preceq \alpha I$ then $\tilde{X}_4 = \beta X_4 \preceq \beta \alpha I = \tilde{\alpha} I$.
\paragraph{Dual feasibility.} Given that $Z \succeq 0$, and by the Schur complement, we have $Z_1 - Z_2^T Z_4^{-1} Z_2 \succeq 0 ,\quad Z_4 \succ 0$. It follows in the same way that those inequalities can be scaled by a positive scalar and, therefore, are equivalent to $
Z_4/\beta \succ 0 ,\quad \beta Z_1  - Z_2^T (Z_4/\beta)^{-1} Z_2 \succeq 0$. By the Schur complement, this means we also have $\tilde{Z} \succ 0$ when $\tilde{Z}$ is given in~\eqref{eq:scaled_XZ}. Next, since $\tilde{Z}_2 = Z_2$, we immediately have $\Vert P^c (\tilde{Z}_2) \Vert_\infty \leq \gamma $.
\paragraph{Zero gradient of the Lagrangian.} If $X_1 = (S-Z_1)^{-1}$, then we scale both sides by $1/\beta$ and obtain $X_1/\beta = ( \beta(S-Z_1))^{-1}$, which is the same as $\tilde{X}_1 = (\tilde{S}-\tilde{Z_1})^{-1}$.
\paragraph{Complementary slackness.} By a simple algebra, we easily confirm that if $ZX =0$ then for any $\beta \neq 0$, we also have $\tilde{X}\tilde{Z} = 0$ where $(\tilde{X},\tilde{Z})$ are given in~\eqref{eq:scaled_XZ}. Moreover, if $Z_4 (X_4 - \alpha I) = 0$ then $(Z_4/\beta) ( \beta X_4 - \beta \alpha I) = \tilde{Z}_4 (\tilde{X}_4 - \tilde{\alpha} I ) = 0$.

In addition to above results, we can check from~\eqref{eq:scaled_XZ} that if $(X,Z)$ has low rank properties: $X_4 = X_1 - X_2^T X_4^{-1}Z_2$ and  $Z_4 = Z_1 - Z_2^T Z_4^{-1}Z_2$, then $(\tilde{X},\tilde{Z})$ also has a low rank. 


\section{Algorithm description}
\label{sec:alg_description}
In this section, we describe the update rules of proximal algorithms used to solve the sparse SEM~\eqref{eq:seml1_primal} when it is arranged in the format of~\eqref{eq:consensus}. ADMM and PPXA algorithms require proximal steps as explained in section~\ref{sec:alg}. An initial solution is selected to be $X_0 = \begin{bmatrix} S^{-1} & 0 \\ 0 & \alpha I \end{bmatrix}$. The main variable is $X$ where $X$ and $X^{+}$ denote the current and next iteration variable, respectively.
\subsection{ADMM algorithm}
 Firstly, we initialize $X = (X_0,X_0,X_0), Z = X_0$ and $Y = (0,0,0) \in \symm^{2n} \times \symm^{2n} \times \symm^{2n}$. Repeat the following steps:
\begin{eqnarray*}
	 X^+_i &=&  \prox_{f_i/\rho} ( Z - Y_i/\rho), \quad \text{for}\;\; i =1,2,3, \\
	 Z^+ & = & (1/3)\sum_{i=1}^3  X^+_i, \\
	Y^+_i &=& Y_i+ \rho(X^+_i - Z^+), \quad \text{for}\;\; i =1,2,3
\end{eqnarray*}
until the primal and dual residual norms are less than a threshold~\cite[\S 7]{BPCP+:10}. 

\subsection{PPXA algorithm}
The update rule follows directly from~\cite[\S 10]{CoP:11} when we choose $\omega_i = 1/3$ for $i=1,2$ and $3$ (uniform.) Firstly, we initialize
$Y=(X_0,X_0,X_0), X=X_0, P = (X_0,X_0,X_0), \bar{P} = (1/3)\sum_{i=1}^3 P_i $. Repeat the following steps:
\begin{eqnarray*}
	 P^+_i &=&  \prox_{\gamma f_i/ \omega_i } ( Y_i), \quad \text{for}\;\; i =1,2,3, \\
	 \bar{P}^+ & = & (1/3)\sum_{i=1}^3  P^+_i, \\
	Y^+_i &=& Y_i+ \lambda( 2\bar{P}^+ - X - P_i), \quad \text{for}\;\; i =1,2,3, \\
X^+ &= & X + \lambda(\bar{P}^+ - X)
\end{eqnarray*}
until a stopping criterion (relative change of cost objective) is less than a threshold.

\section{Brain Region of Interest (ROI)}
\label{app:roi}
\begin{table}[!ht]
\footnotesize

\centering
\caption{Names of region of interests (ROIs) in fMRI connectivity modeling according to Automated Anatomical Labeling (AAL) template.}
\begin{tabular}{|| c | p{6cm} || c | p{6cm} |}   
\hline
No. & Name & No. & name \\ [0.5ex] 
 \hline\hline
1 & Left precentral gyrus (PreCG.L) & 2 &  Right precentral gyrus (PreCG.R)\\
3 & Left superior frontal gyrus (SFGdor.L) & 4 & Right superior frontal gyrus (SFGdor.R)\\
5 & Left superior frontal gyrus, orbital part (ORBsup.L) & 6 & Right superior frontal gyrus, orbital part (ORBsup.R)\\  
7 & Left middle frontal gyrus (MFG.L) & 8 & Right middle frontal gyrus (MFG.R)\\
9 & Left middle frontal gyrus, orbital part (ORBmid.L) & 10 & Right middle frontal gyrus, orbital part (ORBmid.R)\\
11 & Left inferior frontal gyrus, pars opercularis (IFGoperc.L) & 12 & Right inferior frontal gyrus, pars opercularis (IFGoperc.R)\\
13 & Left inferior frontal gyrus, pars triangularis (IFGtriang.L) & 14 & Right inferior frontal gyrus, pars triangularis (IFGtriang.R)\\
15 & Left inferior frontal gyrus, pars orbitalis (ORBinf.L) & 16 & Right inferior frontal gyrus, pars orbitalis (ORBinf.R)\\
17 & Left Rolandic operculum (ROL.L) & 18 & Right Rolandic operculum (ROL.R)\\
19 & Left supplementary motor area (SMA.L) & 20 & Right supplementary motor area (SMA.R)\\
21 & Left olfactory cortex (OLF.L) & 22 & Right olfactory cortex (OLF.R)\\
23 & Left medial frontal gyrus (SFGmed.L) & 24 & Right medial frontal gyrus (SFGmed.R)\\
25 & Left medial orbitofrontal cortex (ORBsupmed.L) & 26 & Right medial orbitofrontal cortex (ORBsupmed.R)\\
27 & Left gyrus rectus (REC.L) & 28 & Right gyrus rectus (REC.R)\\
29 & Left insula (INS.L) & 30 & Right insula (INS.R)\\
31 & Left anterior cingulate gyrus (ACG.L) & 32 & Right anterior cingulate gyrus (ACG.R)\\
33 & Left midcingulate area (DCG.L) & 34 & Right midcingulate area (DCG.R)\\
35 & Left posterior cingulate gyrus (PCG.L) & 36 & Right posterior cingulate gyrus (PCG.R)\\
37 & Left hippocampus (HIP.L) & 38 & Right hippocampus (HIP.R)\\
39 & Left parahippocampal gyrus (PHG.L) & 40 & Right parahippocampal gyrus (PHG.R)\\
41 & Left amygdala (AMYG.L) & 42 & Right amygdala (AMYG.R)\\
43 & Left calcarine sulcus (CAL.L) & 44 & Right calcarine sulcus (CAL.R)\\
45 & Left cuneus (CUN.L) & 46 & Right cuneus (CUN.L)\\
47 & Left lingual gyrus (LING.L) & 48 & Right lingual gyrus (LING.R)\\
49 & Left superior occipital (SOG.L) & 50 & Right superior occipital (SOG.R)\\
51 & Left middle occipital gyrus (MOG.L) & 52 & Right middle occipital gyrus (MOG.R)\\
53 & Left inferior occipital cortex (IOG.L) & 54 & Right inferior occipital cortex (IOG.R)\\
55 & Left fusiform gyrus (FFG.L) & 56 & Right fusiform gyrus (FFG.R)\\
57 & Left postcentral gyrus (PoCG.L) & 58 & Rightpostcentral gyrus (PoCG.R)\\
59 & Left superior parietal lobule (SPG.L) & 60 & Right superior parietal lobule (SPG.R)\\
61 & Left inferior parietal lobule (IPL.L) & 62 &  Right inferior parietal lobule (IPL.R)\\
63 & Left supramarginal gyrus (SMG.L) & 64 & Right  supramarginal gyrus (SMG.R)\\
65 & Left angular gyrus (ANG.L) & 66 & Right angular gyrus (ANG.R)\\
67 & Left precuneus (PCUN.L) & 68 & Right precuneus (PCUN.R)\\
69 & Left paracentral lobule (PCL.L) & 70 & Right paracentral lobule (PCL.R)\\
71 & Left caudate nucleus (CAU.L) & 72 & Right caudate nucleus (CAU.R)\\
73 & Left putamen (PUT.L) & 74 & Right putamen (PUT.R)\\
75 & Left globus pallidus (PAL.L) & 76 & Right globus pallidus (PAL.R)\\
77 & Left thalamus (THA.L) & 78 & Right thalamus (THA.R)\\
79 & Left transverse temporal gyrus (HES.L) & 80 & Right transverse temporal gyrus (HES.R)\\
81 & Left superior temporal gyrus (STG.L) & 82 & Right superior temporal gyrus (STG.R)\\
83 & Left superior temporal pole (TPOsup.L) & 84 & Right superior temporal pole (TPOsup.R)\\
85 & Left middle temporal gyrus (MTG.L) & 86 & Right right middle temporal gyrus (MTG.R)\\
87 & Left middle temporal pole (TPOmid.L) & 88 & Right middle temporal pole (TPOmid.R)\\
89 & Left inferior temporal gyrus (ITG.L) & 90 & Right inferior temporal gyrus (ITG.R)\\
 \hline
\end{tabular}
\label{tab:roi_name}
\end{table}

\end{document}

%% file: zdefs.tex
\newcommand{\BEAS}{\begin{eqnarray*}}
\newcommand{\EEAS}{\end{eqnarray*}}
\newcommand{\BEA}{\begin{eqnarray}}
\newcommand{\EEA}{\end{eqnarray}}
\newcommand{\BEQ}{\begin{equation}}
\newcommand{\EEQ}{\end{equation}}
\newcommand{\BIT}{\begin{itemize}}
\newcommand{\EIT}{\end{itemize}}

\newcommand{\ie}{{\it i.e.}}

\newcommand{\reals}{{\mbox{\bf R}}}
\newcommand{\symm}{{\mbox{\bf S}}}

\newcommand{\Rank}{\mathop{\bf rank}}

\newcommand{\Tr}{\mathop{\bf tr}}
\newcommand{\diag}{\mathop{\bf diag}}

\newcommand{\argmin}{\mathop{\rm argmin}}

\newcommand{\sign}{\mathrm{\bf sign}}

\newcommand{\prox}{\mbox{\bf prox}}

\newcommand{\Cov}{\mathop{\bf cov{}}}
\newcommand{\minimize}{\mathop{\rm minimize{}}}